\documentclass{elsarticle}

\usepackage{amsmath, latexsym}
\usepackage{amssymb}
\usepackage[mathscr]{eucal}
\usepackage{afterpage}
\usepackage{graphicx}
\usepackage{pstricks}
\usepackage{ifthen}
\usepackage{tabularx}
\usepackage{algorithmic}
\usepackage{caption}
\usepackage{subcaption}
\usepackage{textcomp}
\usepackage{algorithm}
\usepackage{siunitx}

\numberwithin{equation}{section}

\graphicspath{{images/}}

\newlength{\halffigwidth}
\newlength{\smallfigwidth}
\newlength{\figheight}
\newcommand{\zero}{\boldsymbol{0}}

\renewcommand{\d}{\boldsymbol{d}}
\newcommand{\e}{\boldsymbol{e}}

\newcommand{\n}{\boldsymbol{n}}

\renewcommand{\u}{\widehat{u}}

\newcommand{\x}{\boldsymbol{x}}
\newcommand{\y}{\boldsymbol{y}}

\newcommand{\m}{\boldsymbol{m}}

\newcommand{\CL}{\mathcal{L}}

\newcommand{\bbC}{\mathbb{C}}
\newcommand{\bbE}{\mathbb{E}}

\newcommand{\bbR}{\mathbb{R}}

\newcommand{\bbV}{\mathbb{V}}

\newcommand{\CB}{\mathcal{B}}

\renewcommand{\to}{\rightarrow}

\renewcommand{\max}{\mathrm{max}}

\usepackage{mathtools}

\DeclarePairedDelimiterX\Basics[1](){ #1}

\newproof{proof}{Proof}

\newcommand{\inc}{\mathrm{inc}}

\definecolor{marker}{rgb}{1,1,0.5}
\definecolor{marker}{rgb}{1,1,0.5}
\definecolor{mydarkgreen}{rgb}{0,0.4,0}
\definecolor{paleyellow}{rgb}{1,1,0.9}

\newcommand{\bsigma}{\boldsymbol{\sigma}}

\newcommand{\revX}[1]{\textcolor{black}{#1}}

\newcommand{\bxi}{\boldsymbol{\xi}}
\newcommand{\bfeta}{\boldsymbol{\eta}}
\newcommand{\btheta}{\boldsymbol{\theta}}
\newcommand{\bzeta}{\boldsymbol{\zeta}}

\journal{Journal of Computational and Applied Mathematics}

\begin{document}

\begin{frontmatter}


\title{An efficient epistemic uncertainty quantification algorithm for a class of stochastic models: A post-processing and domain decomposition framework
}

\author[mg,at]{M. Ganesh}
\ead{mganesh@mines.edu}

\author[sch]{S. C. Hawkins}
\ead{stuart.hawkins@mq.edu.au}

\author[at]{A. M.  Tartakovsky}
\ead{alexandre.tartakovsky@pnnl.gov}

\author[at]{R. Tipireddy\corref{mycorrespondingauthor}}
\cortext[mycorrespondingauthor]{Corresponding author}
\ead{Ramakrishna.Tipireddy@pnnl.gov}

\address[mg]{Department of Applied Mathematics and Statistics, Colorado School of Mines, Golden, CO 80401}
\address[sch]{Department of Mathematics \revX{and Statistics}, Macquarie University, Sydney, NSW 2109, Australia}
\address[at]{Pacific Northwest National Laboratory, Richland, WA 99352}

\date{}
%

\begin{abstract}
Partial differential equations (PDEs) are fundamental for theoretically describing numerous
physical processes that are based on some input fields in   spatial configurations. Understanding the   physical process,
in general, requires computational modeling of the PDE. Uncertainty in the computational model
manifests through lack of precise knowledge of the input field  or configuration. Uncertainty
quantification (UQ)  in the output physical process  is
typically carried out by modeling the uncertainty using a random field, governed
by an appropriate covariance function. This leads to solving   high-dimensional 
stochastic counterparts of the PDE computational models.
Such UQ-PDE models
require a large number of  simulations of the PDE in conjunction with samples in the high-dimensional probability 
space, with probability distribution associated with the covariance function. 
Those UQ computational models having explicit knowledge  of the covariance function
are known as aleatoric  UQ (AUQ) models. The lack of such explicit knowledge leads to epistemic
UQ (EUQ)  models, which typically require solution of a large number of AUQ models. In this article,
using a surrogate, post-processing, and domain decomposition  framework with  coarse stochastic solution 
adaptation, we develop an offline/online algorithm for efficiently simulating a class of EUQ-PDE models.

\end{abstract}


\end{frontmatter}


\section{Introduction}
\label{sec:introduction}

In this work we consider an efficient domain-decomposition-based method, with  coarse stochastic  solution 
adaptation,  for simulation  of a class of stochastic models of the form
\begin{align}
  \begin{split}
  \label{eq:pde}
  \CL v_a(\x,\omega) & = f(\x, \omega),
  \qquad \x \in E, \ \omega \in \Omega,\\
  \CB v_a(\x,\omega) & = g(\x, \omega),
  \qquad \x \in \partial E, \ \omega \in \Omega,\\
  \end{split}
\end{align}
in a (bounded or unbounded) spatial configuration  $E \subseteq \bbR^n$,
where $\CL$ is a partial differential operator and $\CB$ operates  on functions
that are defined on the boundary $\partial E$ of $E$. The stochasticity in the model
manifests through a random field $a$ that may appear as a coefficient of $\CL$,
or may describe the uncertain nature of the configuration $E$. 

Here $\Omega$ is a high-dimensional sample space, and  the dependence of the random field $a: E \times \Omega \to \bbR$
in the stochastic system need not be linear. The  source function $f$ and  boundary data $g$ in~\eqref{eq:pde} are known
data in the model. Deterministic counterparts of the class  of partial differential equations (PDEs) in~\eqref{eq:pde} 
describe numerous physical processes, and several computational models of the determinstic PDE have been widely investigated.

In this work we consider $a(\x,\omega)$
that are normal random fields
with covariance given by
\begin{equation}
  \label{eq:covariance-coefficient}
  C(\x,\y) = \sigma^2
  \exp \left( - | L (\x - \y) |^2  \right),
  \qquad \x,\y \in E,
\end{equation}
for some $n \times n$ matrix $L$ that governs the spatial correlation.
Here $\sigma$ is the standard deviation of $a(\x,\omega)$. 
We emphasize that the common log-normal random field can be considered
a special case of this model in which the random field is incorporated in the
stochastic model as $\exp(a(\x,\omega))$.

Often the  quantity of interest (QoI) is not the solution
$v_a$ of the PDE~\eqref{eq:pde} but some quantity derived from $v_a$.
In some applications the quantity of interest is a functional,
whilst in others the quantity of interest is itself a function of some spatial
variable.
It is therefore convenient to describe the quantity of interest by
$u: D \times \Omega \to \bbC$
where $D \subseteq \bbR^m$ is some appropriate spatial domain.
The regions $D$ and $E$ may or may not coincide, depending on
the application. For example in wave propagation applications, the far-field QoI is a
function of observed direction and hence $D$ is the set of unit vectors.
However the
QoI is obtained from the near-field, which is the solution of the PDE in  a region $E \subseteq \bbR^n$, 
exterior/interior to scattering objects.
We refer to~\cite{nedlec:book,colton:inverse,2019-ghv,2020-dgs,2020-gm}
and references therein for classical and recent literature on forward and inverse acoustic and electromagnetic
wave propagation deterministic and stochastic models.

The standard forward uncertainty  quantification problem, modeled by
the stochastic partial differential equation system~\eqref{eq:pde}, is typically based on
the assumption that the quantities describing the covariance in~\eqref{eq:covariance-coefficient}
are known, leading to  the  {\em aleatoric} UQ (AUQ) problem. 
However, in practice,  sufficient data for precise estimation of $\sigma$ is not available and, hence,
$\sigma$ should be treated as an uncertain parameter, leading to  the associated {\em epistemic} UQ (EUQ) problem.
Over the last two decades, the AUQ-PDE problem has been widely investigated using the
Monte Carlo (MC), quasi-MC (QMC), and generalized polynomial chaos (gPC) 
techniques, see for example~\cite{frances:2, Mai_Knio:stochastic} (and references therein) for the MC, QMC,
and gPC literature for forward AUQ computational models. 
In contrast, the literature on  the EUQ-PDE forward problem
is limited, see the recent work~\cite{epi-2019} and  related EUQ references therein. 

In these published EUQ algorithms, the QoI is assumed to be a scalar and the algorithms were developed accordingly.
In this article, we are interested in QoIs, such as the
far-field, that are functions of spatial variables in $D$. The proposed offline/online approach in this article, with domain-decomposition and coarse stochastic solution framework, 
is entirely different from those considered in the limited computational EUQ literature for the class of stochastic
models described by~\eqref{eq:pde}--\eqref{eq:covariance-coefficient}. Our approach in this article is motivated
by  the gPC stochastic PDE modeling tools developed in our earlier articles~\cite{tipireddy2014basis, gh:uq-fast-multiple, half-plane, tipireddy2017basis, pnnl}.

Each  fixed choice/realization  of the parameter $\sigma$  in the EUQ model leads to one AUQ problem.
The AUQ problem itself is a high-dimensional model in the sampling space $\Omega \subset \bbR^d$, where
the $d$ is typically determined by the decay in the eigenvalues of the covariance.
Since quantifying uncertainties of QoI in the AUQ-PDE  itself is computationally
challenging, the  EUQ problem may even be considered to be computationally
infeasible using the
standard sampling algorithm for the variance parameter  $\sigma \in [\sigma_{\min}, \sigma_{\max}]$ and solving the AUQ-PDE problem for each sample.

The main focus of this article is on developing an efficient algorithm for the EUQ problem. In particular, our approach
may be considered as an offline/online framework, where in the offline part we solve only one AUQ problem with
$\sigma = \sigma_{\max}$,  and using the resulting solution we develop a fast (online) approach to evaluate
the EUQ problem for any large number of samples $\sigma \in [\sigma_{\min}, \sigma_{\max})$, without the
need to further solve the stochastic PDE system. In particular, in addition to quickly obtaining statistical moments for any $\sigma \in [\sigma_{\min}, \sigma_{\max}]$, 
our approach helps to efficiently visualize the QoI for the
EUQ problem through histogram and probability density estimation plots.   The latter can be achieved  using millions of MC samples in $\Omega \subseteq \bbR^d$ for the QoI,
and evaluation of the QoI for any value of $\sigma \in [\sigma_{\min}, \sigma_{\max}]$ with computational cost  essentially
determined by the cost of the single AUQ problem.

The rest of this article is organized as follows.
In the next section, we briefly recall an $N$-term sparse grid gPC (sg-gPC) representation of $u$ for a $d$-dimensional
affine approximation to the random field $a(\x,\omega)$.
In Section~\ref{sec:epistemic} we use the sg-gPC approximation
for $\sigma = \sigma_{\max}$ as a surrogate
to obtain an efficient fast  (online) evaluation algorithm for the EUQ problem.
High-order accuracy of the sp-gPC approximation requires the sparse grid level to depend on $N$ and hence, 
unlike low-order MC/QMC methods, the standard sg-gPC approach requires relatively low stochastic dimension
$d$. In Section~\ref{sec:dd} we recall a recently proposed hybrid of spatial domain decomposition and 
sg-gPC (dd-sg-gPC) for the stochastic dimension reduction. Using a high-order dd-sg-gPC approximation
as the offline surrogate, in Section~\ref{sec:epistemic-dd}, we propose an epistemic dd-sg-gPC algorithm
for the EUQ-PDE stochastic model. In Section~\ref{sec:numerics}, we demonstrate the two epistemic 
algorithms by applying them for EUQ-PDE problems arising in 
a certain class of wave propagation and diffusion models.

\section{gPC approximation}
\label{sec:gpc}

In this section we briefly review how to approximate the solution
of~\eqref{eq:pde}
using a gPC expansion.
The first step is to find a finite-dimensional approximation to
the random coefficient
$a(\x,\omega)$ using a truncated Karhunen-Loeve expansion
\begin{equation}
  \label{eq:kl-coefficient}
  a(\x,\omega) \approx a_0(\x) + \sum_{i=1}^d \sqrt{\lambda_i} a_i(\x)
  \xi_i(\omega)
\end{equation}
where $a_0(\x)$ is the mean of $a(\x)$ and 
$\xi_1,\dots,\xi_d$ are independent Gaussian random variables
with zero mean and unit standard deviation.
The
eigenpairs $(a_i,\lambda_i)$ satisfy the eigenvalue problem
\begin{equation}
  \label{eq:eval-coefficient}
  \int_E C(\x,\y) a_i(\y) \; ds(\y) = \lambda_i a_i(\x),
  \qquad i=1,\dots,d, \ \x \in E.
\end{equation}
In practice~\eqref{eq:eval-coefficient} can be discretized, leading to
an algebraic eigenvalue problem that can be solved efficiently using
the QR algorithm.
All of the eigenvalues of~\eqref{eq:eval-coefficient}
are positive, and we take
$\lambda_1,\dots,\lambda_d$ to be
the $d$ largest eigenvalues, ordered
so that
\begin{displaymath}
\lambda_1 \geq \dots \geq \lambda_d > 0.
\end{displaymath}

Using the truncated expansion~\eqref{eq:kl-coefficient} the
random coefficient $a(\x,\omega)$ is approximated
by a function of the vector valued random variable
$\bxi = (\xi_1,\dots,\xi_d)^T$.
The corresponding gPC approximation to the quantity of interest $u(\x,\bxi)$
is
\begin{equation}
  \label{eq:gpc}
  u_N(\x,\bxi) = \sum_{|\n|=0}^{N} u_{\n}(\x) \psi_{\n}(\bxi)
\end{equation}
where $N$ is the maximum degree of the gPC polynomials and
$u_{\n}(\x)$ are the gPC coefficients, given by
\begin{equation}
  \label{eq:gpc-coefficients}
  u_{\n}(\x) = \langle u(\x,\cdot), \psi_{\n} \rangle.
\end{equation}
Here $\langle \cdot,\cdot \rangle$ is the inner product
\begin{equation}
  \label{eq:inner-product}
  \langle f,g \rangle = \bbE[f g]
  = \int_{\bbR^d} f(\bxi) \overline{g(\bxi)} \; w(\bxi) \; ds(\bxi),
\end{equation}
induced by the Gaussian probability measure
\begin{equation}
  \label{eq:w}
  w(\bxi) = \frac{1}{(2 \pi)^{d/2}} e^{-|\bxi|^2/2}.
\end{equation}
The polynomial basis in~\eqref{eq:gpc}
comprises tensor product polynomials
\begin{equation}
  \label{eq:poly}
\psi_{\n}(\bxi) = \psi_{n_1}(\xi_1) \cdots \psi_{n_d}(\xi_d)
\end{equation}
where $\psi_n$ is the normalized Hermite polynomial of degree $n$,
$\n=(n_1,\dots,n_d)$ is a multi-index 
and $|\n| = n_1 + \dots + n_d$ is the total degree of the tensor product
polynomial.
The tensor product polynomials~\eqref{eq:poly} are orthonormal with
respect to the inner product~\eqref{eq:inner-product}.

The gPC approximation $u_N(\x,\bxi)$ is computationally cheap to evaluate
for any given $\x$ and $\bxi$ because it involves only evaluation of
polynomial terms. (In contrast, direct evaluation of $u(\x,\bxi)$
requires numerical solution of the PDE~\eqref{eq:pde}, which is typically
computationally expensive.)
We will show in Section~\ref{sec:epistemic}
that $u_N(\x,\bxi)$ is a useful surrogate for
$u(\x,\bxi)$ for investigating changes in the solution with respect to the
standard deviation parameter $\sigma$.
The mean and variance of the gPC polynomial,
\begin{equation}
  \label{eq:mean-full}
  \bbE[u_N(\x,\cdot)] = u_{\zero},
  \qquad
  \bbV[u_N(\x,\cdot)] = \sum_{|\n|=1}^{N} | u_{\n}(x) |^2,
\end{equation}
also provide computationally
cheap approximations to the mean and
varaiance of $u(\x,\bxi)$.

In practice we compute the gPC coefficients by approximating the inner product
in~\eqref{eq:gpc-coefficients} using a sparse grid quadrature rule
\begin{equation}
  \label{eq:gpc-quadrature}
  \int_{\bbR^d} f(\bxi) \; w(\bxi) \; ds(\bxi)
  \approx
  \sum_{q=1}^{Q_d} w^{d,Q_d}_q f(\bsigma^{d,Q_d}_q)
\end{equation}
where $\bsigma^{d,Q_d}_q$ and $w^{d,Q_d}_q$
for $q=1,\dots,Q_d$ are the quadrature points and weights respectively.
In our experiments we use a sparse grid quadrature rule
based on a Gauss-Hermite rule
with the number of points $Q_d$ chosen so that the sparse grid level
$\ell = N+2$.
For brevity we suppress the dependence of $Q_d$ on $N$ in our notation.
It follows from~\eqref{eq:gpc-coefficients} and~\eqref{eq:gpc-quadrature}
that assembling the gPC approximation
requires evaluation of $u(\x,\bsigma_q^{d,Q_d})$
for $q=1,\dots Q_d$ by solving the PDE~\eqref{eq:pde}.

\section{Epistemic uncertainty}
\label{sec:epistemic}

In this section we consider the dependence of the quantity of interest
$u(\x,\omega)$ on the standard deviation $\sigma$ of the random field
$a(\x,\omega)$.
For this study, it is convenient to parametrize  the standard
deviation of the random field as $\sigma = \tau \sigma_{\max}$, where $\sigma_{\max}$ is the fixed (and known) maximum value of the standard deviation and
$0 < \tau < 1$ is the epistemic parameter.  

Henceforth, we use 
$C(\x,\y)$ to denote the covariance function
for the  fixed $\sigma_{\max}$,
 and  use $C_{\tau}(\x,\y)$ to represent the epistemic uncertainty in the covariance function.
 Under this parametrization and notation,  
the resulting normal random field $a_\tau$ has covariance
\begin{equation}
  \label{eq:covariance-coefficient-lam}
  C_\tau(\x,\y) = (\tau \sigma_{\max})^2
  \exp \left( - | L (\x - \y) |^2  \right) = \tau^2 C(\x,\y).
\end{equation}

It is convenient to denote the quantity of interest by $u^\tau(\x,\omega)$
where the superscript indicates the dependence 
on the parameter $\tau$.
As in Section~\ref{sec:gpc} we wish to find an approximation to
the mean and standard deviation of 
$u^\tau(\x,\omega)$.
The first step is, again, to find a finite-dimensional approximation to the
random coefficient $a_\tau(\x,\omega)$.
Using the Karhunen-Loeve expansion~\eqref{eq:kl-coefficient} of
$C(\x,\y)$ with
and noting from~\eqref{eq:eval-coefficient} that the eigenvalues of
$C_\tau$ are $\tau^2 \lambda_i$,
we have
\begin{align}
  a_\tau (\x,\omega) & \approx a_0(\x) + 
  \sum_{i=1}^d \sqrt{\tau^2 \lambda_i} a_i(\x) \xi_i(\omega) \nonumber \\
  & \approx a_0(\x) +
  \sum_{i=1}^d \sqrt{\lambda_i} a_i(\x) \tau \xi_i(\omega) \nonumber \\  
  & \approx a_0(\x) +
  \sum_{i=1}^d \sqrt{\lambda_i} a_i(\x) \zeta_i(\omega) \label{eq:kl-coefficient-lambda}
\end{align}
where $\zeta_i = \tau \xi_i$ for $i=1,\dots,d$
are independent Gaussian random variables with zero mean and
standard deviation $\tau$.

Using the truncated expansion~\eqref{eq:kl-coefficient-lambda}
the random coefficient $a_\tau(\x,\omega)$ is approximated by a function
of the vector valued random variable
$\bzeta = (\zeta_1,\dots,\zeta_d)^T = \tau \bxi$.
Similar to Section~\ref{sec:gpc}, we seek an approximation to the quantity
of interest $u^\tau(\x,\bzeta)$ of the form
\begin{equation}
  \label{eq:gpc-lambda}
u^\tau_N(\x,\bzeta) = \sum_{|\n|=0}^{N} u_{\n}^{\tau}(\x) \psi_{\n}^{\tau}(\bzeta)
\end{equation}
where $u^\tau_{\n}(\x)$ are the gPC coefficients, given by
\begin{equation}
  \label{eq:gpc-coefficients-lambda}
  u^\tau_{\n}(\x) = \langle u^\tau(\x,\cdot), \psi^\tau_{\n} \rangle_\tau.
\end{equation}
Here $\langle \cdot,\cdot \rangle_\tau$ is
the inner product
\begin{equation}
  \label{eq:inner-product-lambda}
  \langle f,g \rangle_\tau = \bbE_\tau[f g]
  = \int_{\bbR^d} f(\bzeta) \overline{g(\bzeta)} \; w_\tau(\bzeta) \; ds(\bzeta)
\end{equation}
induced by the Gaussian probability measure
\begin{equation}
  \label{eq:w-lambda}
  w_{\tau}(\bzeta) = \frac{1}{(2 \pi \tau^2)^{d/2}} e^{-\frac{|\bzeta|^2}{2\tau^2}} = \frac{1}{\tau^d} w(\bzeta/\tau).
\end{equation}
The polynomial basis in~\eqref{eq:gpc-lambda} comprises tensor product
polynomials
\begin{displaymath}
  \psi_{\n}^{\tau}(\bzeta) = \psi_{\n}(\bzeta/\tau),
\end{displaymath}
which are easily shown to be orthogonal with respect to the inner
product~\eqref{eq:inner-product-lambda}.

In practice we compute the gPC coefficients
by approximating the inner product in~\eqref{eq:gpc-coefficients-lambda}
using a sparse grid quadrature rule with points
$\tau \bsigma^{d,Q_d}_q$ and weights $w^{d,Q_d}_q$ for
$q = 1,\dots,Q_d$.
This requires evaluation of $u^\tau(\x,\tau \bsigma^{d,Q_d}_q)$
for $q=1,\dots,Q_d$ by solving the PDE~\eqref{eq:pde}.

\subsection*{Surrogate gPC approximation}

In the remainder of this section we suppose that the
gPC approximation $u_N(\x,\bxi)$ described in
Section~\ref{sec:gpc} has been computed.
From~\eqref{eq:kl-coefficient-lambda} we have
\begin{displaymath}
a_\tau(\x,\bzeta) = a(\x,\bzeta)
\end{displaymath}
so that
it is appropriate to use $u_N(\x,\bzeta)$ as a surrogate for
$u^\tau(\x,\bzeta)$.
Using this surrogate
in~\eqref{eq:gpc-coefficients-lambda} gives
surrogate gPC coefficients
\begin{equation}
  \label{eq:gpc-coefficients-surrogate}
  \u^\tau_{\n}(\x)
  = \langle u_N(\x,\cdot), \psi^\tau_{\n} \rangle_\tau.
\end{equation}
Using the expansion~\eqref{eq:gpc} of $u_N(\x,\cdot)$
in~\eqref{eq:gpc-coefficients-surrogate} we derive
\begin{equation}
  \label{eq:times-T}
  \u^{\tau}_{\n}(\x) = \sum_{|\m| = 1}^N T^{d,\tau}_{\n,\m} u_{\m}(\x),
\end{equation}
where
\begin{align}
  T^{d,\tau}_{\n,\m}
  & = \langle \psi_{\m},\psi^{\tau}_{\n} \rangle_\tau \nonumber \\
  & = \int_{\bbR^d} \psi_{\m}(\bzeta) \psi^{\tau}_{\n}(\bzeta) \; w_\tau(\bzeta) \; ds(\bzeta)\nonumber \\
  & = \int_{\bbR^d} \psi_{\m}(\bzeta) \psi_{\n}(\bzeta/\tau) \; w_\tau(\bzeta) \; ds(\bzeta)\nonumber \\
  & = \int_{\bbR^d} \psi_{\m}(\tau \bxi) \psi_{\n}(\bxi) \; w(\bxi) \; ds(\bxi). \label{eq:T}
\end{align}
In the last line we have used the change of variables
$\bzeta = \tau \bxi$.
In practice the integral in~\eqref{eq:T} is approximated using the
sparse grid quadrature rule~\eqref{eq:gpc-quadrature}.

Using the details above,
our surrogate gPC approximation to $u^\tau(\x,\bzeta)$ is
\begin{equation}
  \label{eq:gpc-epistemic}
  \u^{\tau}_N(\x,\bzeta)
  = \sum_{|\n| = 0}^N \u^{\tau}_{\n}(\x) \psi^\tau_{\n}(\bzeta)
\end{equation}
and the surrogate coefficients $\u^{\tau}_{\n}(\x)$ are computed from the
coefficients $u_{\n}(\x)$ using the matrix-vector product~\eqref{eq:times-T}.
Thus the surrogate gPC coefficients are computationally inexpensive to compute
and, in particular, the surrogate gPC coefficients can be obtained without
computing further solutions of the PDE~\eqref{eq:pde}.
The mean and the variance of the surrogate gPC polynomial,
\begin{equation}
  \label{eq:mean-gpc-surrogate}
\bbE_{\tau}[\u^{\tau}_N(\x,\cdot)] = \u^{\tau}_{\zero},
  \qquad
\bbV_{\tau}[\u^{\tau}_N(\x,\cdot)] = \sum_{|\n|=1}^{N} | \u_{\n}^{\tau}(\x) |^2
\end{equation}
also provide computationally cheap approximations 
to the mean
and variance
$u^\tau(\x,\bzeta)$.

\section{Domain decomposition method with coarse stochastic solution}
\label{sec:dd}

\begin{figure}
  \centering
  \includegraphics[width=10cm]{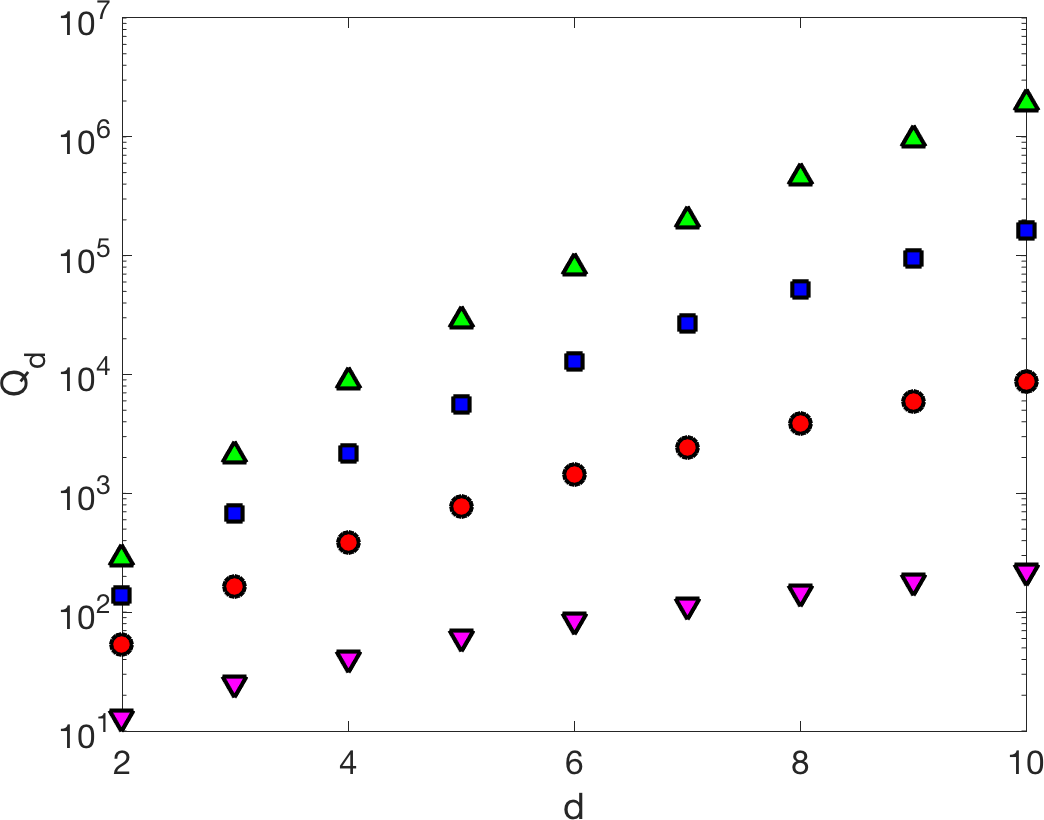}
  \caption{\label{fig:sg-points}
    Plot showing the number of sparse grid points $Q_d$ against dimension
    $d$ for $N=1$ ($\triangledown$),
    $N=3$ ($\bigcirc$), $N=5$ ($\square$) and $N=7$ ($\triangle$).}
\end{figure}

The accuracy of the gPC representation in Section~\ref{sec:gpc}
depends on taking $d$ sufficiently large that $a(\x,\omega)$ is well
approximated by the truncated Karhunen-Loeve
expansion~\eqref{eq:kl-coefficient}. On the other hand, reducing the computational cost of  the gPC method requires keeping $d$ sufficiently small such that the number of sparse grid
quadrature points, needed to evaluate the coefficients
in~\eqref{eq:gpc-coefficients}, is manageable.
In Figure~\ref{fig:sg-points} we demonstrate the rapid growth in the number
of sparse grid points required to compute the gPC approximation~\eqref{eq:gpc}
as the stochastic dimension $d$ increases for polynomials
with maximum degree $N=1,3,5,7$.

In this section we briefly review a recently proposed domain decomposition
method in conjunction with coarse stochastic solution~\cite{pnnl} that provides optimal low dimensional subspaces of
$\bbR^d$ on which to approximate the solution $u(\x,\omega)$.
This facilitates the construction of a new gPC approximation with lower
stochastic dimension than $d$, and hence significantly fewer stochastic
quadrature points are required (see Figure~\ref{fig:sg-points}).

The key to the domain decomposition method is to observe that the
quantity of interest $u(\x,\omega)$
can itself be viewed as a random field on $D$ with covariance
\begin{equation}
  \label{eq:covariance-u}
C_u(\x,\y) = \bbE\bigg[ \Big(u(\x,\cdot) - \bbE[u(\x,\cdot)]\Big)\Big(u(\y,\cdot) - \bbE[u(\y,\cdot)]\Big)\bigg], \qquad \x,\y \in D.
\end{equation}
It follows that the Hilbert-Karhunen-Loeve expansion can be used to find a low
stochastic-dimension approximation to $u(\x,\omega)$, similar to
the expansion of $a(\x,\omega)$ in Section~\ref{sec:introduction}.
Motivated by the fact that random fields typically have lower rank
approximation on smaller
domains~\cite{pnnl},
we decompose the spatial domain $D$ into subdomains $D_1,\dots,D_S$
satisfying
\begin{displaymath}
  D = \bigcup_{s=1}^{S} D_s, \qquad D_s \cap D_{s'} = \emptyset \ 
  \mbox{for $s \neq s'$}.
\end{displaymath}
On each subdomain we use the readily computable approximation
to the covariance in~\eqref{eq:covariance-u}
\begin{equation}
  \label{eq:covariance-u1}
  C^s_{u_1}(\x,\y) = \bbE\bigg[ \Big(u_1(\x,\cdot) - \bbE[u_1(\x,\cdot)]\Big)\Big(u_1(\y,\cdot) - \bbE[u_1(\y,\cdot)]\Big)\bigg], \qquad \x,\y \in D_s,
\end{equation}
where the coarse stochastic solution $u_1(\x,\bxi)$ is the degree 1 gPC approximation given by~\eqref{eq:gpc}.
Using the orthonormality of the tensor product
polynomials~\eqref{eq:poly} and the expansion~\eqref{eq:gpc} we have
\begin{equation}
  \label{eq:covariance-u1-series}
  C^s_{u_1}(\x,\y) = \sum_{|\n| = 1} u_{\n}(\x) u_{\n}(\y),
  \qquad \x,\y \in D_s.
\end{equation}
For polynomial degree $N=1$
the number of sparse grid points required is very low
(see Figure~\ref{fig:sg-points}) so that computing the 
gPC approximation $u_1(\x,\bxi)$ is
computationally inexpensive, even when $\bxi$ has the full
$d$ dimensions.

Using~\eqref{eq:covariance-u1-series}, 
the truncated Karhunen-Loeve expansion for $u_1(\x,\bxi)$ is
\begin{equation}
  \label{eq:kl-u}
  u_1(\x,\bxi(\omega)) \approx u_{\zero}(\x) + \sum_{i=1}^d \sqrt{\mu^s_i} b^s_i(\x)
  \eta^s_i(\omega),
\end{equation}
where $u_{\zero}(\x)$ is the mean of $u_1(\x,\cdot)$
and 
$\eta^s_1,\dots,\eta^s_d$ are independent Gaussian random variables
with zero mean and unit standard deviation.
The eigenpairs $(b^s_i,\mu^s_i)$ satisfy the eigenvalue problem
\begin{equation}
  \label{eq:eval-coefficient-dd}
  \int_D C^s_{u_1}(\x,\y) b^s_i(\y) \; ds(\y) = \mu^s_i b^s_i(\x),
  \qquad i=1,\dots,d.
\end{equation}
In practice,~\eqref{eq:eval-coefficient-dd} can be discretized
(see example-specific details in Section~\ref{sec:numerics}),
leading to an algebraic eigenvalue problem that can be solved efficiently
using the QR algorithm.
All of the eigenvalues of~\eqref{eq:eval-coefficient-dd} are positive,
and we take $\mu_1^s,\dots,\mu_d^s$ to be the $d$ largest eigenvalues,
ordered so that
\begin{displaymath}
\mu^s_1 \geq \dots \geq \mu^s_d > 0.
\end{displaymath}

On the left hand side of~\eqref{eq:kl-u}, the random field $u_1(\x,\cdot)$
is represented
as a function of $\xi_1(\omega),\dots,\xi_d(\omega)$.
On the right hand side of~\eqref{eq:kl-u}, the random field $u_1(\x,\cdot)$ is
represented as a function of $\eta^s_1(\omega),\dots,\eta^s_d(\omega)$.
That is, we have two different parametrizations of $u_1(\x,\cdot)$ with
respect to $d$ stochastic variables.
These parametrizations are related via the orthogonal linear
transformation~\cite{pnnl}
\begin{equation}
  \label{eq:transformation}
  \eta^s_i = \sum_{j=1}^d a^s_{ij} \xi_j,
  \qquad i=1,\dots,d,
\end{equation}
where
\begin{equation}
  \label{eq:aij}
  a^s_{ij} = \frac{1}{\sqrt{\mu^s_i}} \int_{D_s} u_{\e_j}(\x) b^s_i(\x) \; ds(\x),
\end{equation}
and $\e_j$ is the $j$th Euclidean vector in $\bbR^d$.

The approximation~\eqref{eq:kl-u}
is associated with the restricted domain $D_s$ on which 
the variance of $u_1(\x,\cdot)$ is expected to be less than
on the full domain $D$.
Consequently the eigenvalues $\mu^s_1,\dots,\mu^s_d$
are expected to decay sufficiently quickly that the number of terms on the
right hand side of~\eqref{eq:kl-u} can be reduced without significantly
compromising the
accuracy of the approximation.
This motivates further truncating the
Karhunen-Loeve expansion~\eqref{eq:kl-u} to obtain the approximation
\begin{equation}
  \label{eq:kl-u-reduced}
  u_1(\x,\bxi(\omega)) \approx u_{\zero}(\x) + \sum_{i=1}^r \sqrt{\mu^s_i} b^s_i(\x)
  \eta^s_i(\omega),
\end{equation}
with $r < d$.
The expansion~\eqref{eq:kl-u-reduced}
provides an approximation of the random field
$u_1(\x,\cdot)$  for $\x \in D_s$ using a reduced dimension stochastic space.

The advantage of~\eqref{eq:kl-u-reduced} is in
identifying a reduced dimension stochastic
space appropriate for constructing a local approximation to $u(\x,\cdot)$
for $\x \in D_s$.
In particular,
writing $\bfeta^s = (\eta^s_1,\dots,\eta^s_r)$ we have
\begin{equation}
  \label{eq:approx-dd}
  u(\x,\bxi(\omega)) \approx u^s(\x,\bfeta^s(\omega)),
  \qquad \x \in D_s,
\end{equation}
where
$u^s(\x,\bfeta^s) = u(\x,A^s \bfeta^s)$ and
\begin{equation}
  \label{eq:eta-to-xi}
  \bfeta^s = (A^s)^T \bxi.
\end{equation}
Here $A^s$ is the $d \times r$ matrix with entries
\begin{displaymath}
A^s_{ij} = a^s_{ij},\qquad i=1,\dots,d, \ j=1,\dots,r,
\end{displaymath}
where $a^s_{ij}$ is given by~\eqref{eq:aij}.

Following the details in Section~\ref{sec:gpc} we approximate
$u^s(\x,\bfeta^s)$ by a truncated gPC approximation
\begin{equation}
  \label{eq:gpc-dd}
u^s_{N_s}(\x,\bfeta^s) = \sum_{|\n|=0}^{N_s} u^s_{\n}(\x) \psi_{\n}(\bfeta^s)
\end{equation}
where
\begin{equation}
  \label{eq:gpc-coefficients-dd}
  u^s_{\n}(\x) = \langle u^s(\x,\bfeta^s), \psi_{\n} \rangle
  = \langle u(\x,A^s \bfeta^s), \psi_{\n} \rangle.
\end{equation}
In practice the inner product in~\eqref{eq:gpc-coefficients-dd} is
computed using
the sparse grid quadrature rule
\begin{equation}
  \label{eq:gpc-quadrature-dd}
  \int_{\bbR^r} f(\bxi) \; w(\bfeta) \; ds(\bfeta)
  \approx
  \sum_{q=1}^{Q_r} w^{r,Q_r}_q f(\bsigma^{r,Q_r}_q).
\end{equation}
The gPC approximation~\eqref{eq:gpc-dd}
is cheaper to compute than~\eqref{eq:gpc} because it has lower
stochastic dimension and hence the sparse grid quadrature rule has fewer
points (see Figure~\ref{fig:sg-points}).
In particular, $Q_r < Q_d$.

Using the gPC approximation~\eqref{eq:gpc-dd} we compute approximations
to the mean and variance,
\begin{equation}
  \label{eq:mean-dd}
  \bbE[u(\x,\cdot)] \approx u_{\zero}^{s},
  \qquad
  \bbV[u(\x,\cdot)] \approx \sum_{|\n|=1}^{N_s} | u_{\n}^{s}(\x) |^2,
    \qquad \x \in D_s.
\end{equation}

\section{Epistemic uncertainty for the Domain Decomposition method}
\label{sec:epistemic-dd}

As in Section~\ref{sec:epistemic}, we again let the standard deviation
of the random field $a_\tau$
be $\tau \sigma_{\max}$ for some parameter $\tau > 0$.
Then using the Karhunen-Loeve expansion~\eqref{eq:kl-coefficient-lambda}
for $a_\tau$,
the change of variables $\btheta^s = \tau \bfeta^s$,
and the transformation~\eqref{eq:eta-to-xi}, we can approximate
\begin{equation}
  \label{eq:approx-dd-epistemic}
  u^\tau(\x,\bzeta) \approx u^{\tau,s}(\x,\btheta^s), \qquad
  \x \in D_s,    
\end{equation}
where
$u^{\tau,s}(\x,\btheta^s) = u^\tau(\x,A^s \btheta^s)$, with
\begin{equation}
  \label{eq:theta-to-zeta}
  \btheta^s = (A^s)^T \bzeta,
\end{equation}
and $A^s$ is the matrix given in Section~\ref{sec:dd}.

Similar to Section~\ref{sec:gpc}, we seek an approximation to the local
quantity of interest $u^{\tau,s}(\x,\btheta^s)$
for $\x \in D_s$ of the form
\begin{equation}
  \label{eq:gpc-lambda-s}
  u^{\tau,s}_N(\x,\btheta^s) = \sum_{|\n|=0}^{N_s} u^{\tau,s}_{\n}(\x)
  \psi^\tau_{\n}(\btheta^s)
\end{equation}
where $u^{\tau,s}_{\n}(\x)$ are the gPC coefficients, given by
\begin{equation}
  \label{eq:gpc-coefficients-lambda-s}
  u^{\tau,s}_{\n}(\x) = \langle u^{\tau,s}(\x,\cdot), \psi^\tau_{\n} \rangle_\tau.
\end{equation}

\subsection*{Surrogate gPC approximation}

As in Section~\ref{sec:epistemic}, in the remainder of this section
we suppose that the local gPC approximation $u^s(\x,\bfeta)$ has been
computed.
Following the details in Section~\ref{sec:epistemic}, it is appropriate to use
$u^s_N(\x,\btheta^s)$ as a surrogate for $u^{\tau,s}_N(\x,\btheta^s)$.
Using this surrogate in~\eqref{eq:gpc-coefficients-lambda-s} gives surrogate
gPC coefficients
\begin{equation}
  \label{eq:gpc-coefficients-surrogate}
  \u^{\tau,s}_{\n}(\x)
  = \langle u^s_N(\x,\cdot), \psi^\tau_{\n} \rangle_\tau.
\end{equation}
Similar to~\eqref{eq:times-T} in Section~\ref{sec:epistemic} we have
\begin{equation}
  \label{eq:times-T-s}
  \u^{\tau,s}_{\n}(\x) = \sum_{|\m| = 1}^N T^{r,\tau}_{\n,\m} u^s_{\m},
\end{equation}
where $T^{r,\tau}$ is given by~\eqref{eq:T}.

Using the details above, our surrogate approximation to
$u^{\tau,s}(\x,\btheta^s)$ for $\x \in D_s$ is
\begin{equation}
  \label{eq:gpc-epistemic-dd}
  \u^{s,\tau}_N(\x,\btheta^s)
  = \sum_{|\n| = 0}^N \u^{s,\tau}_{\n}(\x) \psi^\tau_{\n}(\btheta^s)
\end{equation}
and the surrogate coefficients $\u^{s,\tau}_{\n}(\x)$ are computed
from the coefficients $u^s_{\n}(\x)$ using the matrix-vector
product~\eqref{eq:times-T-s}.
As in Section~\ref{sec:epistemic}, the surrogate gPC coefficients are
computationally inexpensive to compute and, in particular,
the surrogate coefficients can be obtained without computing further
solutions of the PDE~\eqref{eq:pde}.
The mean and variance of the surrogate gPC polynomial,
\begin{equation}
    \label{eq:mean-dd-surrogate}
\bbE_{\tau}[\u^{s,\tau}_N(\x,\cdot)] = \u^{s,\tau}_{\zero},
  \qquad
\bbV_{\tau}[\u^{s,\tau}_N(\x,\cdot)] = \sum_{|\n|=1}^{N} | \u_{\n}^{s,\tau}(\x) |^2
\end{equation}
also provide computationally cheap approximations to the mean and variance
of $u^{\tau,s}(\x,\btheta^s)$ for $\x \in D_s$.

\section{Numerical Results}
\label{sec:numerics}
In this section we demonstrate the efficiency of our surrogate based EUQ algorithm
for two distinct stochastic models. The first model describes wave propagation in an
unbounded medium, exterior to an uncertain configuration, 
and the second model is a widely investigated diffusion process
with uncertain permeability field. We investigate our EUQ algorithm 
with and without applying the domain decomposition framework, and show that the latter  is
essential for solving a $100$-dimensional model, arising due to slow decay of eigenvalues. 

\subsection*{A stochastic wave propagation model in unbounded region}

\begin{figure}
  \centering
  \includegraphics[width=10cm]{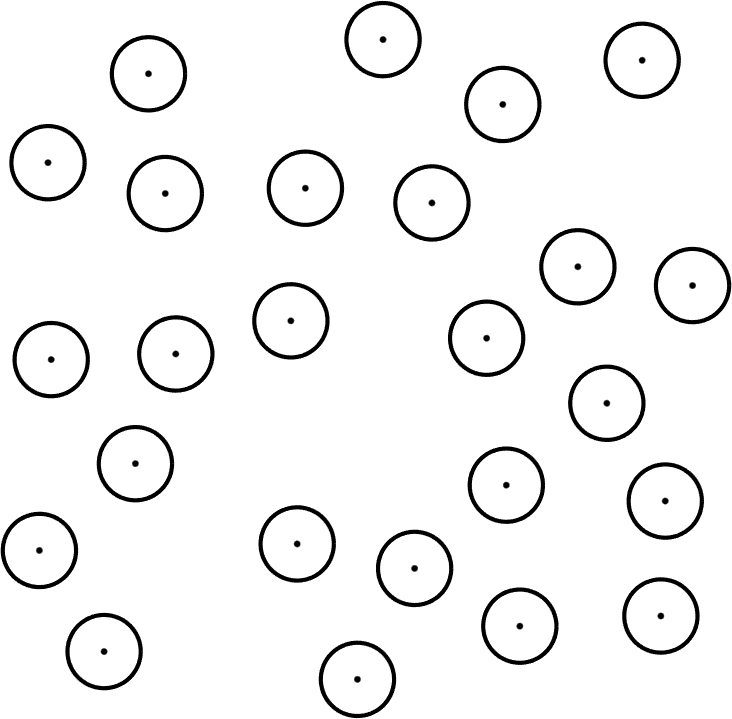}
  \caption{\label{fig:scatterers10}
    Visualization of a sample of $K(\omega)$ for
    a configuration of 25 cylinders.
    }
\end{figure}

We consider multiple scattering of a time harmonic transverse electric
(TE) polarized electromagnetic wave by a configuration of
$M=25$ parallel perfectly conducting
cylinders whose positions are described by a random field. For this uncertain
configuration model, we let $K(\omega) \subseteq \bbR^2$ denote the cross section of the cylinders.
Then
\begin{equation}
  \label{eq:scatterers}
  K(\omega) = \cup_{i=1}^M K_i(\omega),
\end{equation}
where the unit disk $K_i(\omega)$ is the cross section of the $i$th
cylinder.
In our model $K_i(\omega)$ has center $\x_i(\omega) = \x_i + a(\x_i,\omega) \d$,
where $\d = (1,1)^T / \sqrt{2}$ is the  translation direction 
and
$a(\x,\omega)$ is a normal random
field satisfying~\eqref{eq:covariance-coefficient}
on $E=\bbR^2$
with variance $\sigma_{\max} = 0.05$ and
$L = \mathop{\mathrm{diag}}(1/l,1/l)$.
Here, $l=10$ is the correlation length
and the constant $\x_i$ is the mean of the centre $\x_i(\omega)$ of
the $i$th scatterer.
In our experiments $\x_i$ for $i=1,\dots,M$ are chosen
at random in the square
$[0,10] \times [0,10]$.
A sample of $K(\omega)$ is visualized in Figure~\ref{fig:scatterers10}.

\begin{figure}
  \centering
  \includegraphics[width=10cm]{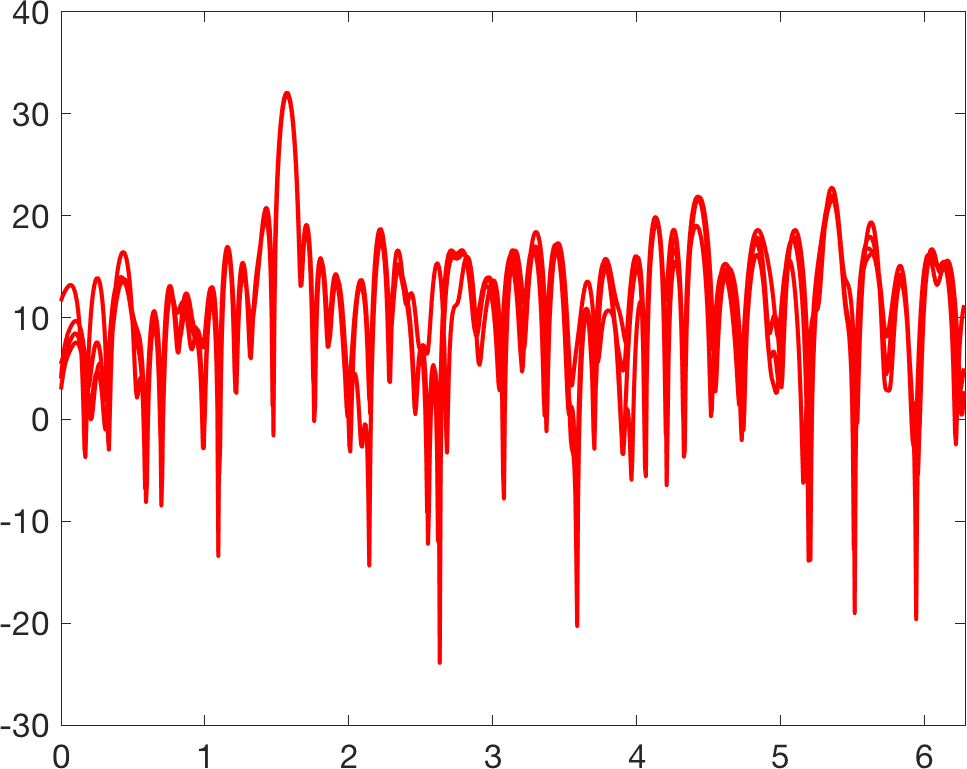}
  \caption{\label{fig:farfields10}
    Visualization of the radar cross section
    of the configuration in Figure~\ref{fig:scatterers10}
    in decibels
    $10 \log_{10} 2 \pi u(\theta;\omega)$.
  }
\end{figure}

The configuration $K(\omega)$ is illuminated by an incident plane wave
modelled by the scalar field
\begin{displaymath}
v^\inc(\x) = e^{ik x_2},
\end{displaymath}
where $k = c/\omega$ is the wavenumber and $c$ is the speed of light.
The resulting scattered field is modelled by scalar field
$v(\x;\omega)$ that satisfies the two-dimensional Helmholtz
equation
\begin{equation}
  \label{eq:helmholtz}
  \triangle v(\x;\omega) + k^2 v(\x;\omega) = 0, \qquad
  \x \in \bbR^2 \setminus \overline{K(\omega)},
\end{equation}
together with the perfect conductor boundary condition
\begin{equation}
  \label{eq:bc}
  v(\x;\omega) = - v^\inc(\x), \qquad \x \in \partial K(\omega),
\end{equation}
where $\partial K(\omega)$ denotes the boundary of $K(\omega)$,
and the radiation condition
\begin{equation}
  \label{eq:radiation-condition}
  \lim_{r \to \infty} \sqrt{r} \left( \frac{\partial v}{\partial r}
  - i k v \right) = 0, \qquad r = |\x|,
\end{equation}
uniformly for all directions $\theta$,
where we have used two dimensional
polar coordinates $\x = (r \cos \theta,r \sin \theta)$.
A consequence of~\eqref{eq:radiation-condition} is that the scattered
field can be written
\begin{equation}
  \label{eq:farfield}
  v(r,\theta,\omega) \approx \frac{e^{ikr}}{\sqrt{r}} v_\infty(\theta;\omega),
  \qquad \theta \in D,
\end{equation}
where $D = [0,2 \pi) \subseteq \bbR$.
  The function $v_\infty(\theta,\omega)$ is known as the far field of $v$
  and is used to compute the radar cross section
  \begin{equation}\label{eq:cross-section-db}
    u(\theta,\omega) = |v_\infty(\theta;\omega)|^2,
  \end{equation}
  which is typically
  the quantity of interest in applications.
  The radar cross section
    of the configuration in Figure~\ref{fig:scatterers10}
  is visualized in Figure~\ref{fig:farfields10}.

\begin{figure}
  \centering
  \includegraphics[width=10cm]{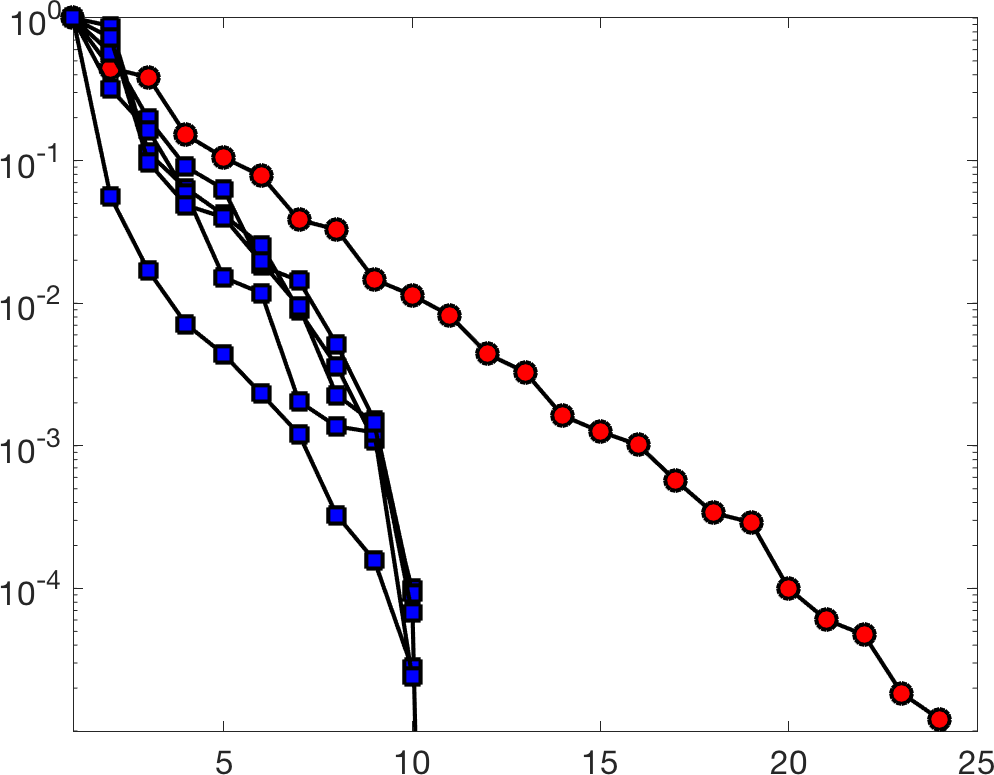}
  \caption{\label{fig:eigenvalue5}
    Decay of the eigenvalues of
    the covariance function $C(\x,\y)$ ($\bigcirc$)
    and the covariance functions $C^s_{u_1}(\x,\y)$ ($\square$) for
    $s=1,\dots,S$ with
    stochastic dimension $d=10$ and
    $S=5$ subdomains.
  }
\end{figure}

  We discretize the eigenvalue problem~\eqref{eq:eval-coefficient} using
  a Monte Carlo quadrature scheme with equal weights and nodes at the
  mean centers $\x_i$ for $i=1,\dots,25$ of the cylinders.
  The decay of the eigenvalues of $C(\x,\y)$ is shown in
  Figure~\ref{fig:eigenvalue5} and we choose truncation parameter $d=10$
  leading to a Karhunen-Loeve approximation~\eqref{eq:kl-coefficient}
  of $a(\x,\omega)$ with respect to the random variable $\bxi$
  with stochastic dimension $d=10$.

  In Table~\ref{table:full-error} we present the relative
  error maximum norm of the 
  approximation~\eqref{eq:mean-full} to the mean computed using the
  gPC scheme in Section~\ref{sec:gpc}.
  The gPC coefficients are computed using~\eqref{eq:gpc-coefficients}
  and the quadrature scheme~\eqref{eq:gpc-quadrature}.
  For each fixed $\bxi$ we compute $u(\x,\bxi)$
  in~\eqref{eq:gpc-coefficients} by
  efficiently solving the associated deterministic
  wave scattering
  PDE~\eqref{eq:helmholtz}--\eqref{eq:radiation-condition} and computing the
  far field~\eqref{eq:farfield}
  using the Matlab package MieSolver~\cite{miesolver}, which is based on the
  Mie series~\cite{rayleigh1881,mie1908} (see also the
  books~\cite{hulst:light,bohren,rother}).
  In practice we compute a discrete approximation to the maximum norm
  using 1000 equally spaced points in $[0,2\pi]$,
  \begin{equation}
    \label{eq:farfield-points}
      x_i = 2\pi \frac{i-1}{1000}, \qquad i=1,\dots,1000.
    \end{equation}

    Next we decompose the domain $D = [0,2\pi]$ into $S$ subdomains
    \begin{displaymath}
      D_s = [2\pi(s-1)/S,2\pi s/S), \qquad s=1,\dots,S.
    \end{displaymath}
    On each domain $D_s$ we compute the gPC approximation $u_1(\x,\bxi)$
    as described above.
    We discretize the eigenvalue problem~\eqref{eq:eval-coefficient-dd} using
    an equal weight quadrature scheme with
    nodes at the points~\eqref{eq:farfield-points}.
    The decay of the eigenvalues of $C^s_{u_1}(\x,\y)$ for $S=5$ subdomains
    is shown in 
    Figure~\ref{fig:eigenvalue5}.
    The figure shows that the eigenvalues of $C^s_{u_1}(\x,\y)$ decay faster
    than the eigenvalues of $C(\x,\y)$ so that it is appropriate to apply
    the dimension reduction approach in Section~\ref{sec:dd}.

\begin{figure}
  \centering
  \begin{tabular}{ccc}
    \includegraphics[width=4cm]{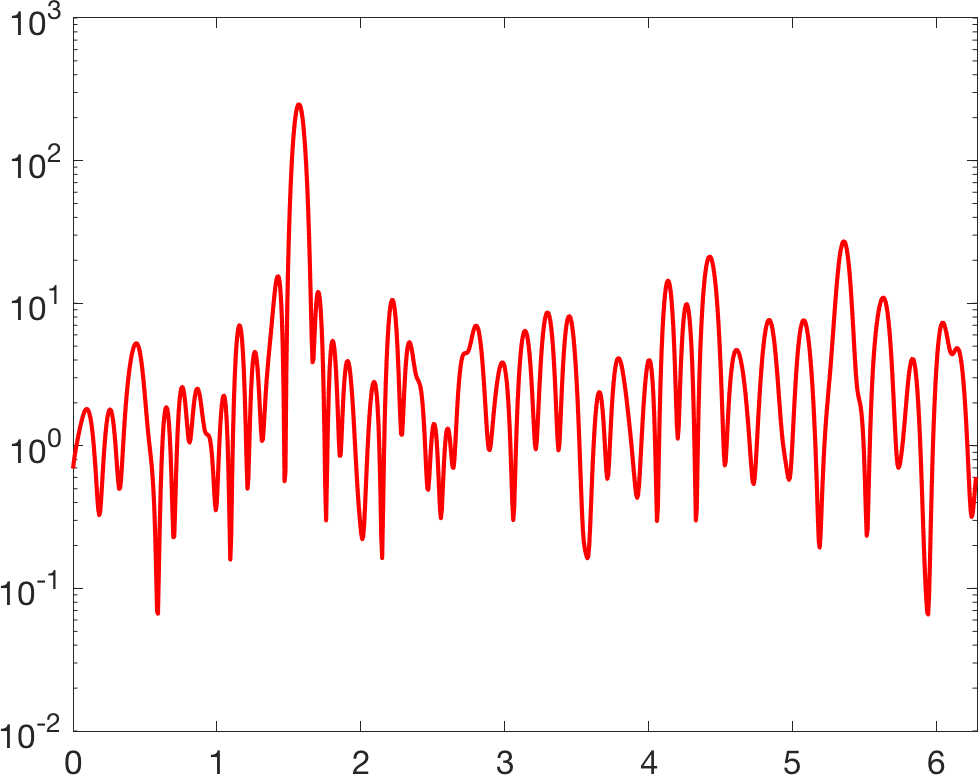}
    &
    \includegraphics[width=4cm]{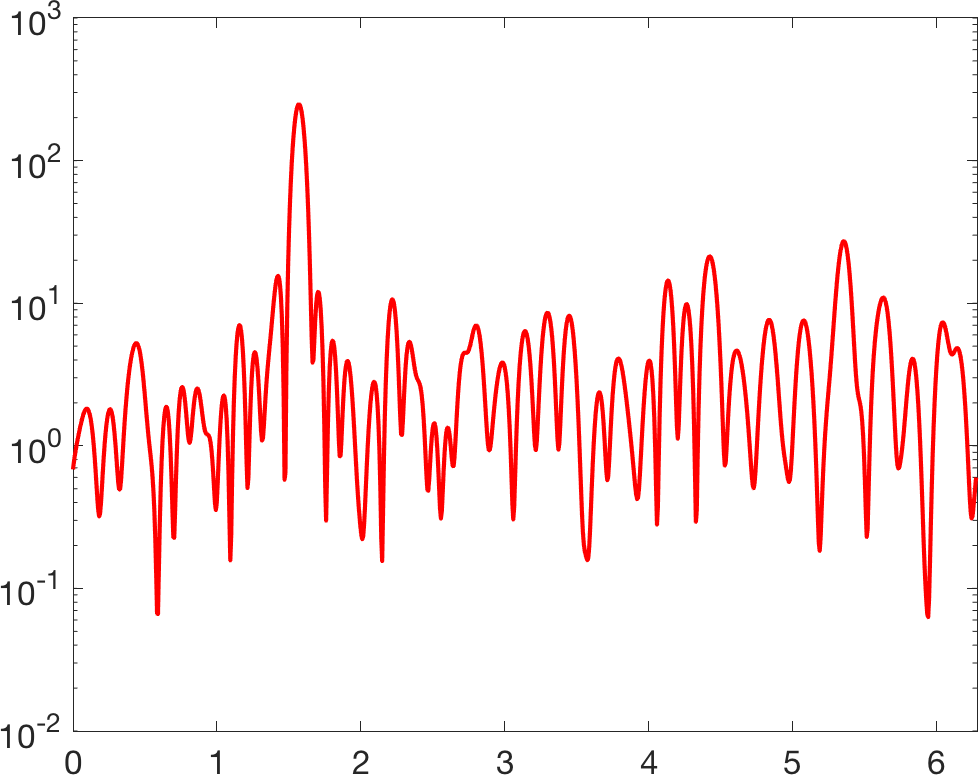}
    &
    \includegraphics[width=4cm]{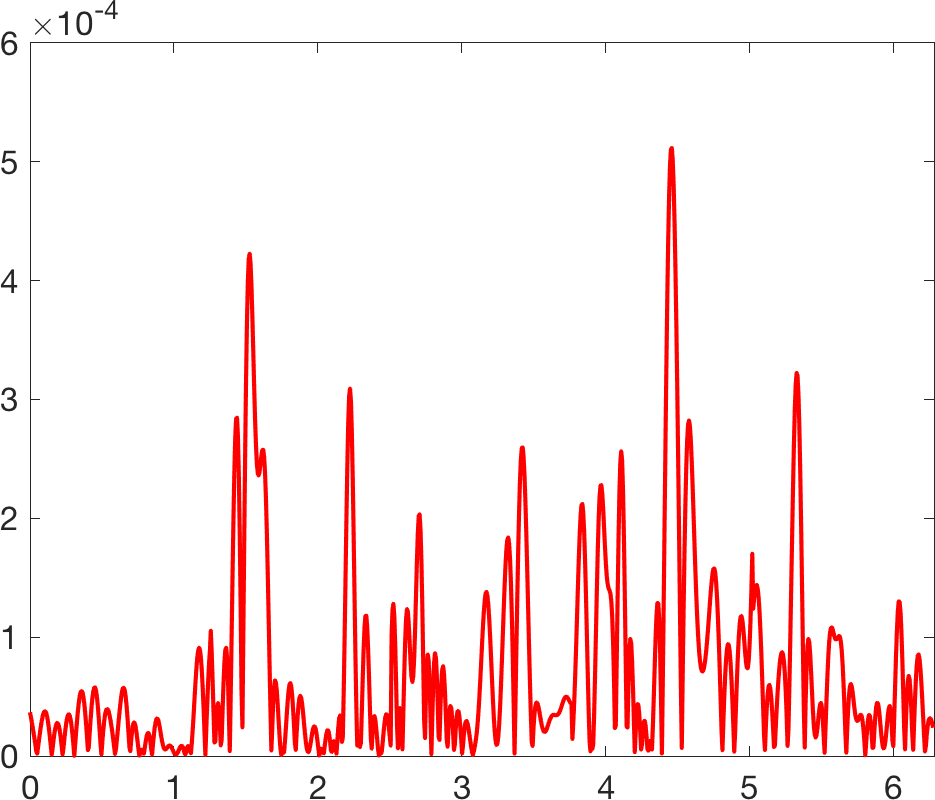}\\
    \parbox{4cm}{\centering
      (a) \\
      Mean of $u_N(\x,\bxi)$ for $d=10$}
    &
    \parbox{4cm}{\centering
      (b)\\
      Mean of $u_N^s(\x,\bzeta)$ for $r=5$}
    &
    \parbox{4cm}{\centering
      (c)\\
      relative difference}
  \end{tabular}
  \caption{\label{fig:mean5}
    Approximations to the mean of $u(\x,\omega)$ obtained
    without domain decomposition (a) and with domain
    decomposition (b) using $S=5$ subdomains.}
\end{figure}

\begin{figure}
  \centering
  \begin{tabular}{ccc}
    \includegraphics[width=4cm]{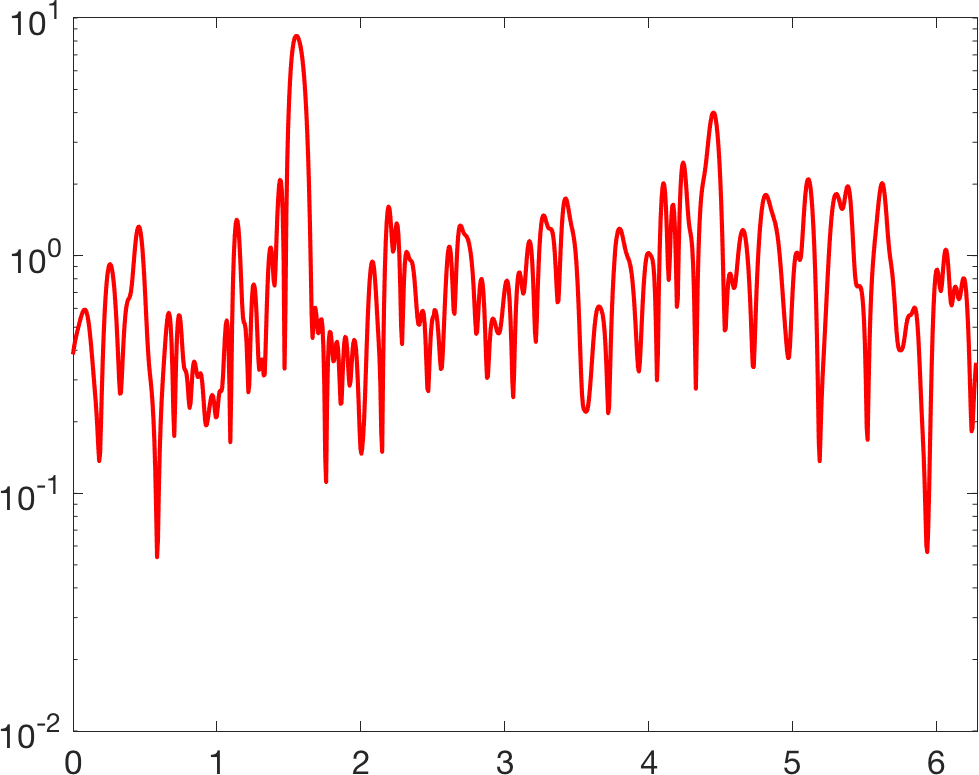}
    &
    \includegraphics[width=4cm]{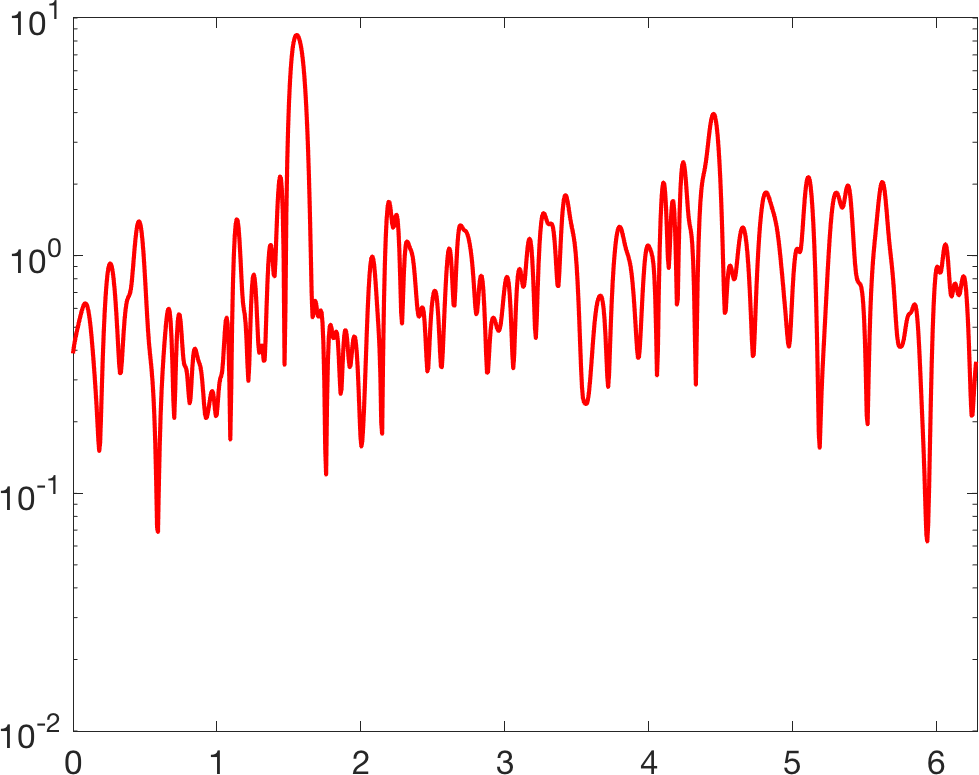}
    &
    \includegraphics[width=4cm]{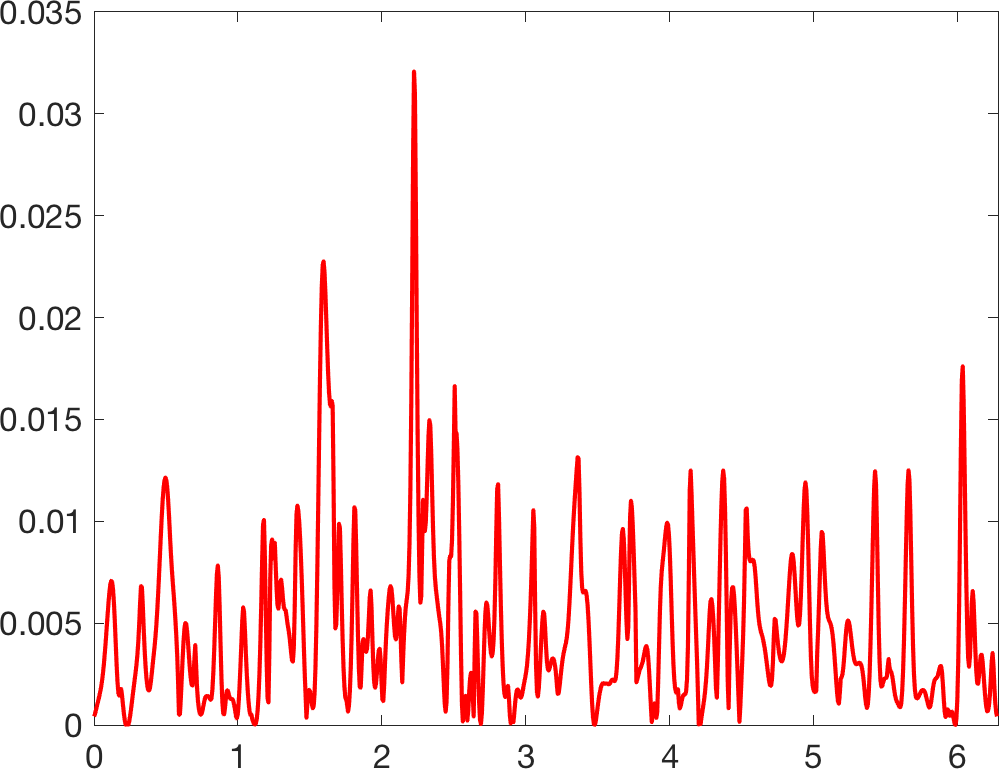}\\
    \parbox{4cm}{\centering
      (a) \\
      Standard deviation of $u_N(\x,\bxi)$ for $d=10$}
    &
    \parbox{4cm}{\centering
      (b)\\
      Standard deviation of $u_N^s(\x,\bzeta)$ for $r=5$}
    &
    \parbox{4cm}{\centering
      (c)\\
      relative difference}
  \end{tabular}
  \caption{\label{fig:var5}
    Approximations to the standard deviation of $u(\x,\omega)$ obtained
    without domain decomposition (a) and with domain
    decomposition (b) using $S=5$ subdomains.}
\end{figure}

  In Table~\ref{table:dd-error} we tabulate the relative
  error maximum norm of the 
  approximation~\eqref{eq:mean-dd} to the mean computed using the
  domain decomposition scheme in Section~\ref{sec:dd}.
  We demonstrate that the domain decomposition technique is able to
  produce an approximation of comparable quality to the $d=10$ gPC
  scheme using only $r=5$ dimensions.
  The $d=10$ gPC scheme requires evaluation of the
  PDE model for 162\,025 values of the stochastic parameter $\bxi$,
  whereas the domain decomposition scheme with $r=5$ requires evaluation of
  the PDE at only 5593 values of the stochastic parameter $\bfeta$
  for each subdomain.
  Thus the domain decomposition method is significantly cheaper.
  In Figures~\ref{fig:mean5} and~\ref{fig:var5} we plot the mean and
  standard deviations of the cross
  sections~\eqref{eq:cross-section-db} computed using
  the gPC method~\eqref{eq:mean-full} with $d=10$ and the
  domain decomposition method~\eqref{eq:mean-dd} with $r=5$.

  In Table~\ref{eq:table-epistemic-dd}
  we demonstrate the accuracy of our approximations by comparing with
  reference solutions obtained using Monte Carlo applied directly for
  random fields $a(\x,\omega)$ with covariance given
  by~\eqref{eq:covariance-coefficient} with $\sigma$ replaced by
  $\tau \sigma_{\max}$. In Tables~\ref{eq:table-wave-gpc-epistemic-err-d10}--\ref{eq:table-wave-gpc-dd-epistemic-err-d10}
  we demonstrate  the accuracy of the epistemic post-processing algorithm  for the wave propagation
  problem, when compared to the associated aleatoric gPC model, with
  and without domain decomposition.

\begin{figure}
  \centering
  \includegraphics[width=10cm]{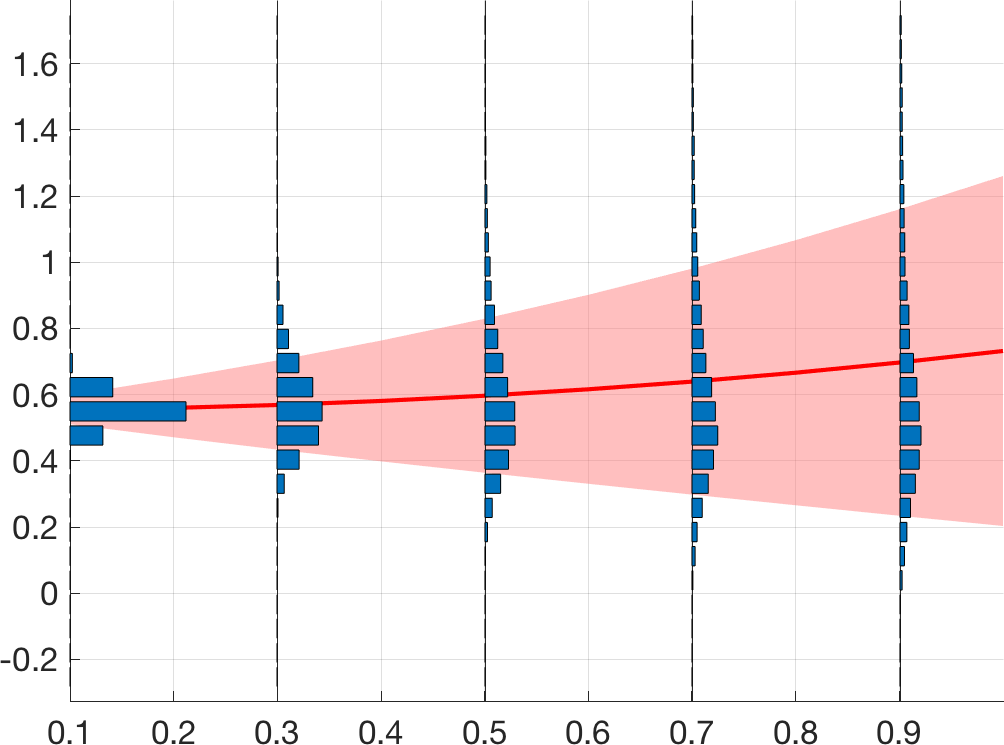}
  \caption{\label{fig:epistemic-dd}
    Visualization of the mean (solid line) and standard deviation
    (indicated by shading) of the cross section $u(\x,\omega)$
    plotted against $\tau$.
    The random
    field $a(\x,\omega)$ has standard
    deviation $\tau \sigma_{\max}$ for
    $0<\tau<1$.
    The mean and standard deviation are computed using the
    efficient method in Section~\ref{sec:epistemic-dd}.
    Samples with 10\,000 Monte Carlo samples are visualized in the histograms.}
\end{figure}

 In Figure~\ref{fig:epistemic-dd}
  we visualize the approximations to the
  mean and standard deviation of the backscattered cross section
  as a function of $\tau$, where $\tau \sigma_{\max}$ is the standard
  deviation of the input random field $a(\x,\omega)$.
  The backscattered cross section approximations are 
  computed using~\eqref{eq:mean-dd-surrogate}
  with backscattering direction $\x = 3 \pi/2$.
  
\begin{table}
  \centering
  \begin{tabular}{cc}
    \hline
    $d$ & rel.~error \\
    \hline
    4 & 2.55e-03 \\
    6 & 1.00e-03 \\
    8 & 3.70e-04 \\
    10 & 1.60e-04 \\
    \hline
  \end{tabular}
  \caption{\label{table:full-error}
    Maximum norm relative error in the approximation to the mean $u(\x,\omega)$
    computed using the gPC method with
    stochastic dimensions $r=4,6,8,10$.
    The reference solution is computed using 1\,048\,576 Monte Carlo simulations.}
\end{table}

\begin{table}
  \centering
  \begin{tabular}{cc}
    \hline
    $d$ & rel.~error \\
    \hline
    4 & 1.20e-03\\
    5 & 6.62e-04\\
    6 & 4.15e-04\\
    \hline
  \end{tabular}
  \caption{\label{table:dd-error}
    Maximum norm error of the approximation to the mean $u(\x,\omega)$
    computed using the domain decomposition method with
    $S=5$ subdomains and stochastic dimensions $r=4,5,6$.
    The reference solution is computed using 1\,048\,576
    Monte Carlo simulations.}
\end{table}

\begin{table}
  \centering
  \begin{tabular}{ccc}
    \hline
    $\tau$ & rel.~error \\
    \hline
    0.2 &	2.22e-03\\
    0.4 &	1.80e-03\\
    0.6 &	1.10e-03\\
    0.8 &	3.59e-04\\
    1.0 &	6.62e-04\\
    \hline
  \end{tabular}
  \caption{\label{eq:table-epistemic-dd}
    Maximum norm relative error of the approximation to the mean computed
    using the domain decomposition method in Section~\ref{sec:epistemic-dd}
    for random field
    with standard deviation $\tau \sigma_{\max}$.
    The reference solution is computed using 1\,048\,576
    Monte Carlo simulations.}
\end{table}

\begin{table}
  \centering
  \begin{tabular}{ccc}
    \hline
    $\tau$ & rel.~error \\
    \hline
0.2 & 	1.71e-02\\
0.4 & 	1.33e-02\\
0.6 & 	8.06e-03\\
0.8 & 	2.81e-03\\
1.0 & 	1.82e-04\\
    \hline
  \end{tabular}
  \caption{\label{eq:table-epistemic-full}
    Maximum norm relative error of the approximation to the mean computed
    using the method in Section~\ref{sec:epistemic}
    for random field
    with standard deviation $\tau \sigma_{\max}$.
    The reference solution is computed using 1\,048\,576
    Monte Carlo simulations.}
\end{table}

\begin{table}[!ht]
 \centering
 \begin{tabular}{cccc}
   \hline
   $\tau$ & max. rel.~mean error & max. rel.~variance error\\
   \hline
   0.9 & 9.9416e-04 & 3.7682e-02 \\
   0.8 & 2.9073e-03 & 1.0801e-01 \\
   0.7 & 5.3900e-03 & 1.9246e-01 \\
   0.6 & 8.1333e-03 & 2.7751e-01 \\
   0.5 & 1.0871e-02 & 3.5289e-01 \\
   0.4 & 1.3383e-02 & 4.1282e-01 \\
   0.3 & 1.5497e-02 & 4.5643e-01 \\
   0.2 & 1.7087e-02 & 4.8519e-01 \\
    0.1 & 1.8072e-02 & 5.0129e-01 \\
   \hline
 \end{tabular}
 \caption{\label{eq:table-wave-gpc-epistemic-err-d10}
   {\bf Ten-dimensional stochastic diffusion model}:
   Maximum norm relative error of the approximate mean and variance computed using the
   gPC approach and the post-processed surrogate method
   in Section~\ref{sec:epistemic}.}
\end{table}

\begin{table}[!ht]
 \centering
 \begin{tabular}{cccc}
   \hline
   $\tau$ & max. rel.~mean error & max. rel.~variance error\\
   \hline
   0.9 & 2.4060e-04 & 6.4024e-03 \\
   0.8 & 5.5901e-04 & 1.7486e-02 \\
   0.7 & 9.1775e-04 & 3.5164e-02 \\
   0.6 & 1.2785e-03 & 5.6954e-02 \\
   0.5 & 1.6084e-03 & 8.0734e-02 \\
   0.4 & 1.8847e-03 & 1.0409e-01 \\
   0.3 & 2.0956e-03 & 1.2477e-01 \\
   0.2 & 2.2395e-03 & 1.4090e-01 \\
0.1 & 2.3216e-03 & 1.5108e-01 \\
   \hline
 \end{tabular}
 \caption{\label{eq:table-wave-gpc-dd-epistemic-err-d10}
   {\bf Ten-dimensional stochastic diffusion model}:
   Maximum norm relative error of the approximate mean and variance computed using the
   domain decomposition approach and the post-processed surrogate method
   in Section~\ref{sec:epistemic-dd}.}
\end{table}

\newpage

\subsection*{A stochastic diffusion model}
We consider the EUQ counterpart of  the  stochastic diffusion example investigated in~\cite{pnnl}.
The large, but bounded, spatial domain of the model is  $E  = [0, 240] \times [0, 60]$. For  $\x=(x_1,x_2) \in E$,
stochasticity in the diffusion model with a mixed (Dirichlet  and Neumann) boundary condition (on vertical
and horizontal boundaries of $E$),  is induced by a log-normal field:
\begin{align}\label{eq:diff_ex}
	-\nabla \cdot \left[\exp(a(\x,\omega)) \nabla u(\x,\omega)\right]&=f(\x) \;\; \rm{in}~E\times \Omega,  \nonumber \\  
	u(\x,\omega)=50 \;\; \text{on}~x_1 = 0, \qquad & \; \; 
	u(\x,\omega)=25 \;\; \text{on}~x_1 = 240, \nonumber \\
	\frac{\partial}{\partial \boldsymbol{n}} u(\x,\omega) = 0  \;\; \text{on}~x_2 = 0, \qquad & 
	\frac{\partial}{\partial \boldsymbol{n}}  u(\x,\omega)  = 0  \;\; \text{on}~x_2 = 60,
\end{align}
where $\boldsymbol{n}(\x)$ denotes the unit outward normal to $E$ at $\x$ and
the discontinuous   function $f$ is such that the sink value of $f$  at the center of the domain $E$ is $-1$, and 
$f$ is  zero elsewhere. For this model, the QoI is the solution $u$ of~\eqref{eq:diff_ex} and hence $D = E$. As in~\cite{pnnl}, we choose 
$L = \mathop{\mathrm{diag}}(1/24, 1/20)$, and  the mean and  maximum standard deviation of $\exp(a)$ to be
5.0 and 2.5 respectively.  Hence the mean and maximum standard deviation
of the normal random field $a$ are respectively
\[
a_0(x) = \ln \left (\frac{5}{\sqrt{1+(\frac{2.5}{5})^2}} \right ) = \ln (5/\sqrt{1.25}); 
\qquad 
\sigma_{\max} = \sqrt{\ln (1.25)}.
\]
We demonstrate the accuracy of the post-processed surrogate solutions of the above model with
standard deviations $\tau \sigma_{\max}$, for $0  < \tau < 1 $, obtained using the gPC 
solution of the model with standard deviation $\sigma_{\max}$.

A $10$-dimensional
truncated version of the stochastic  model was simulated
for fixed standard deviation $\sigma$
using
the gPC approximation with  $N = 3$ and the sparse grid level $N+2$
in~\cite{pnnl}.
Simulation of reference solutions for this problem requires solving the
diffusion PDE $8761$ times, once for each sparse-grid point in
the ten-dimensional stochastic space. Using this gPC
approximation as the reference solution, it was shown in~\cite{pnnl} that a domain 
decomposition (DD) version of the solution has accuracy $6 \times 10^{-5}$ in mean and 
$7 \times 10^{-3}$ in variance. Such an accurate gPC-DD solution was obtained  
with $S = 8$ sub-domains and, for   $s =1, \dots, S$,  in each sub-domain $D_s$ the stochastic dimension was
chosen to be $r = 3$ and the gPC $N_s = N$. Below, we use the same parameters 
for the EUQ model.

We note that it takes only a few seconds of the CPU time to post-process
the $\tau=1$  case reference solution $u_N^\tau$ to compute the approximation    
$\widehat{u}_N^\tau$ to the solution   $u_N^\tau$, for any $0 < \tau < 1$.
In particular, this post-processing does not requires further solves of the
PDE model.

\begin{figure}[]
  \centering
  \begin{tabular}{ccc}
    \includegraphics[width=4cm]{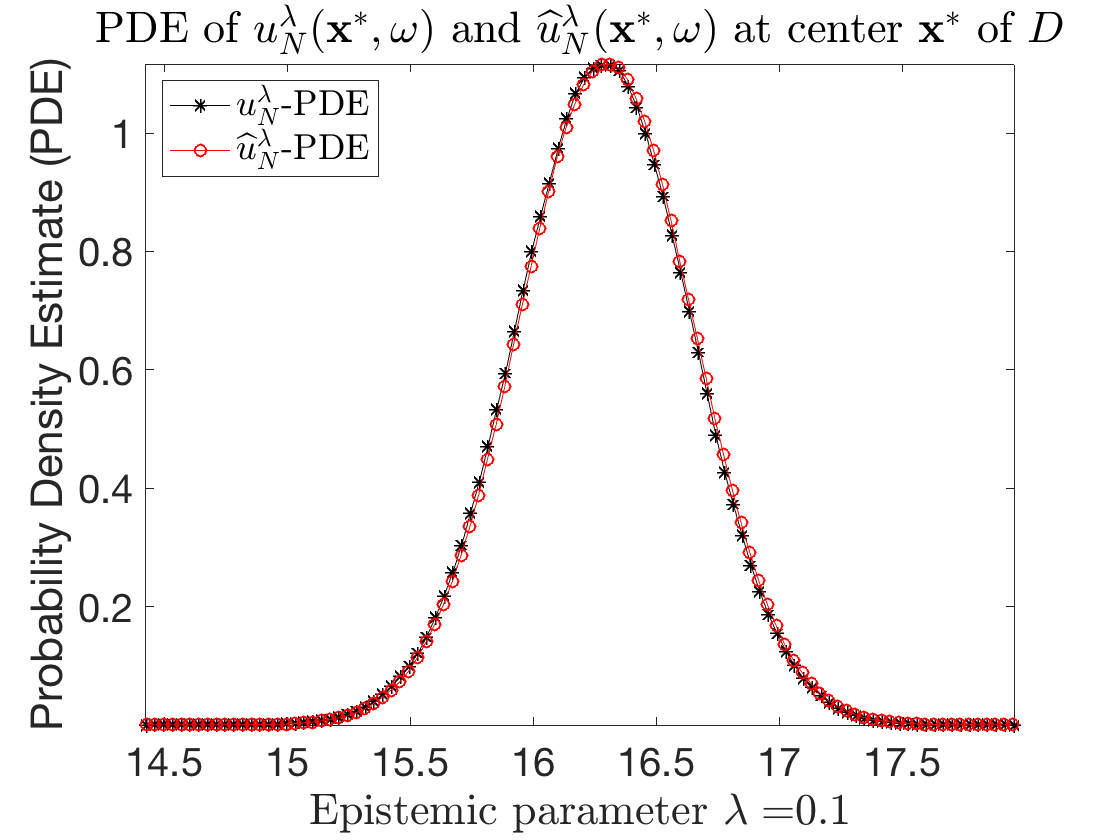}
    &
    \includegraphics[width=4cm]{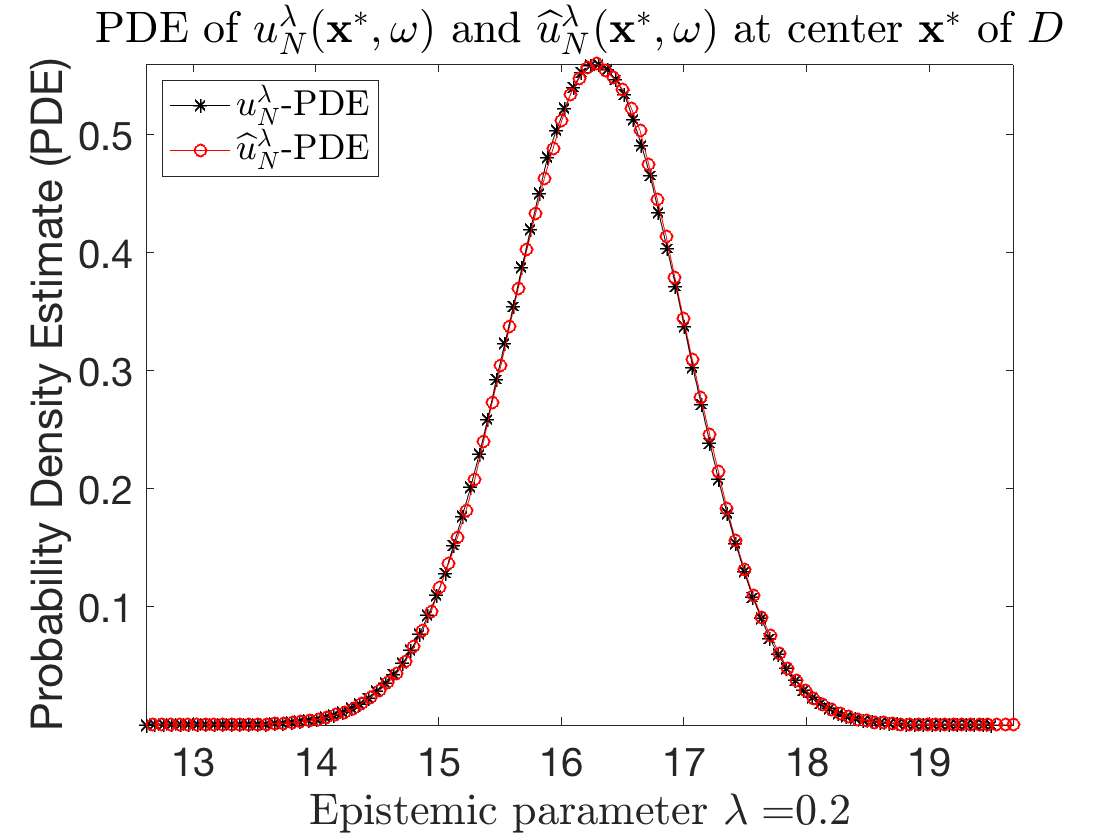}
    &
    \includegraphics[width=4cm]{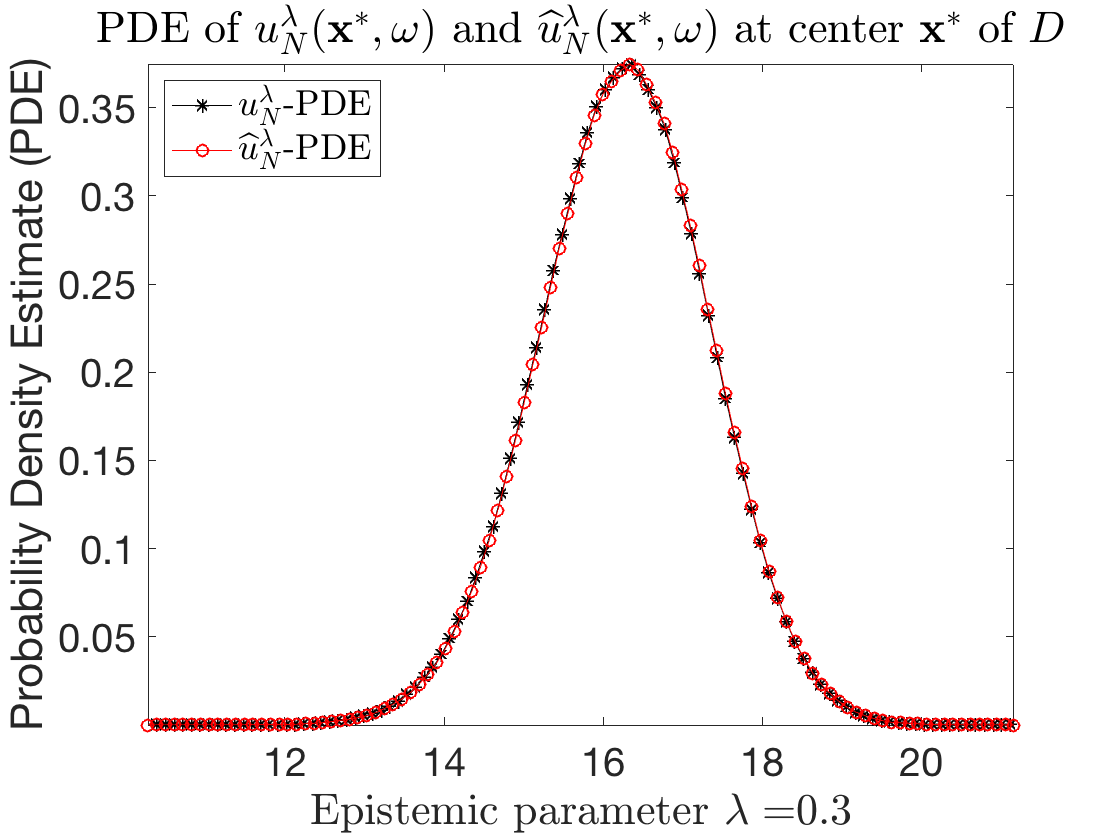}\\
    \parbox{4cm}{\centering
      (a) $\tau = 0.1$}
    &
    \parbox{4cm}{\centering
      (b) $\tau = 0.2$}
    &
    \parbox{4cm}{\centering
      (c) $\tau = 0.3$} \\
      \includegraphics[width=4cm]{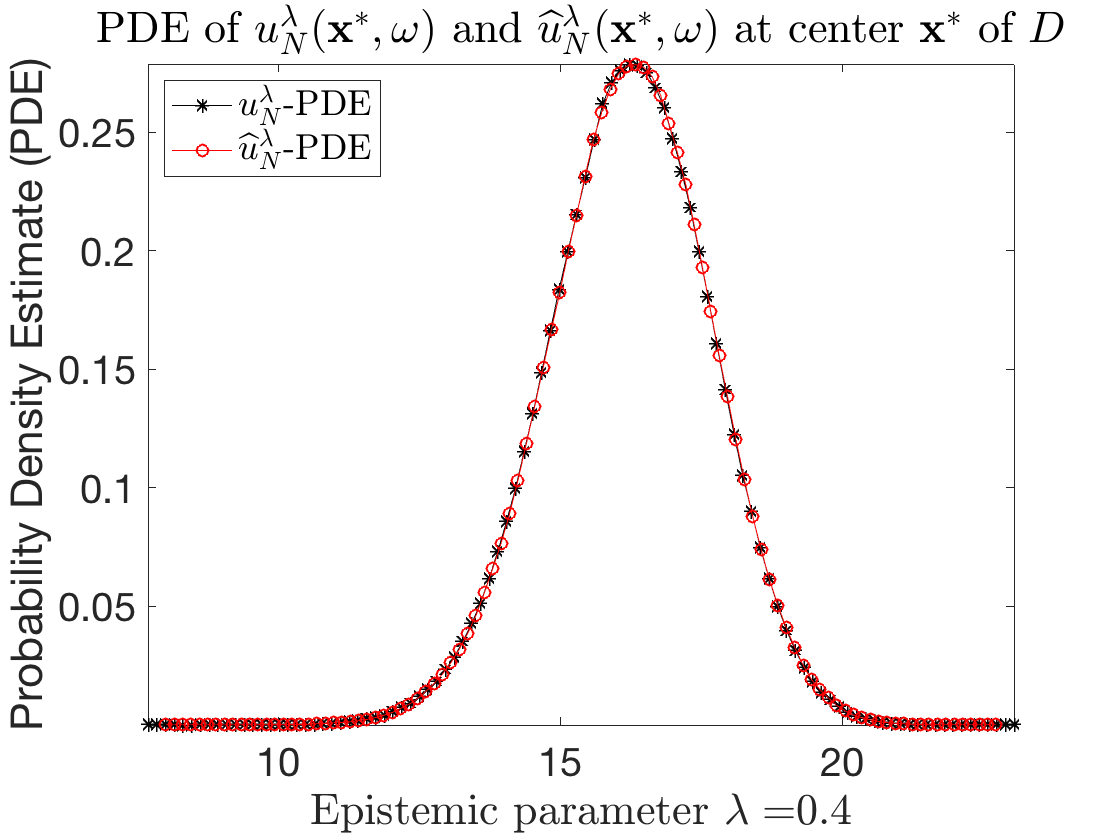}
    &
    \includegraphics[width=4cm]{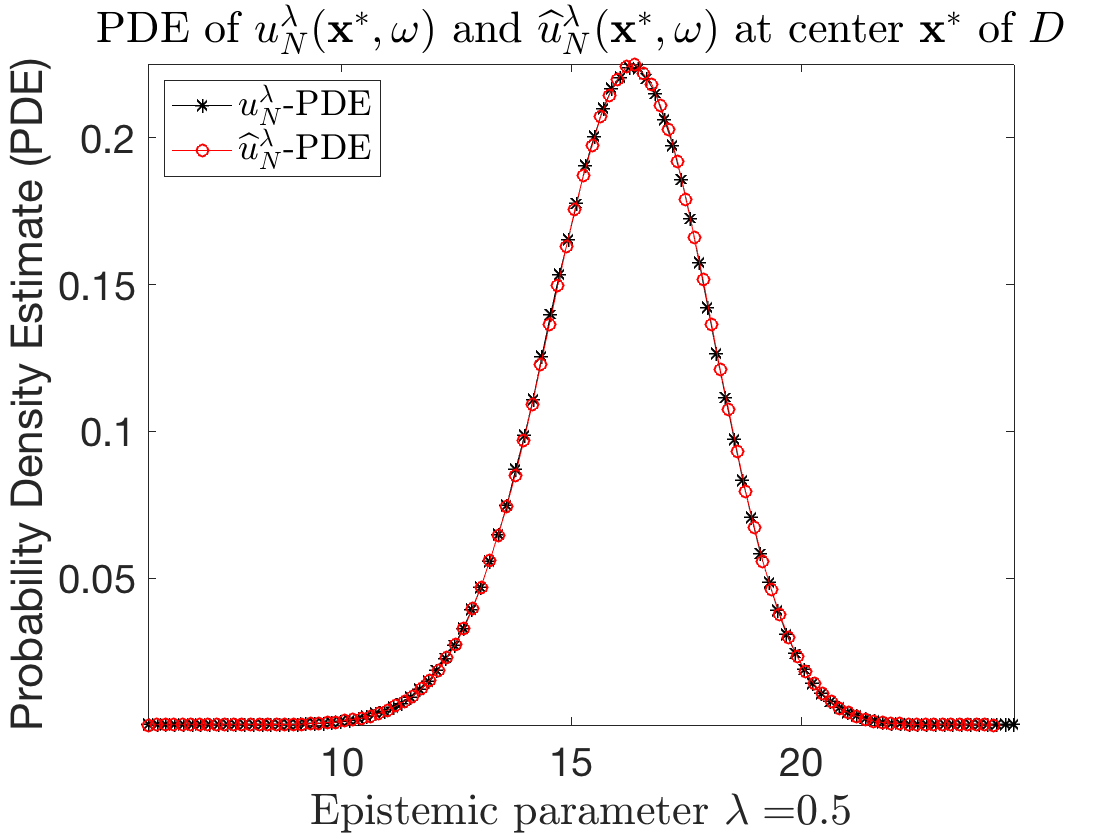}
    &
    \includegraphics[width=4cm]{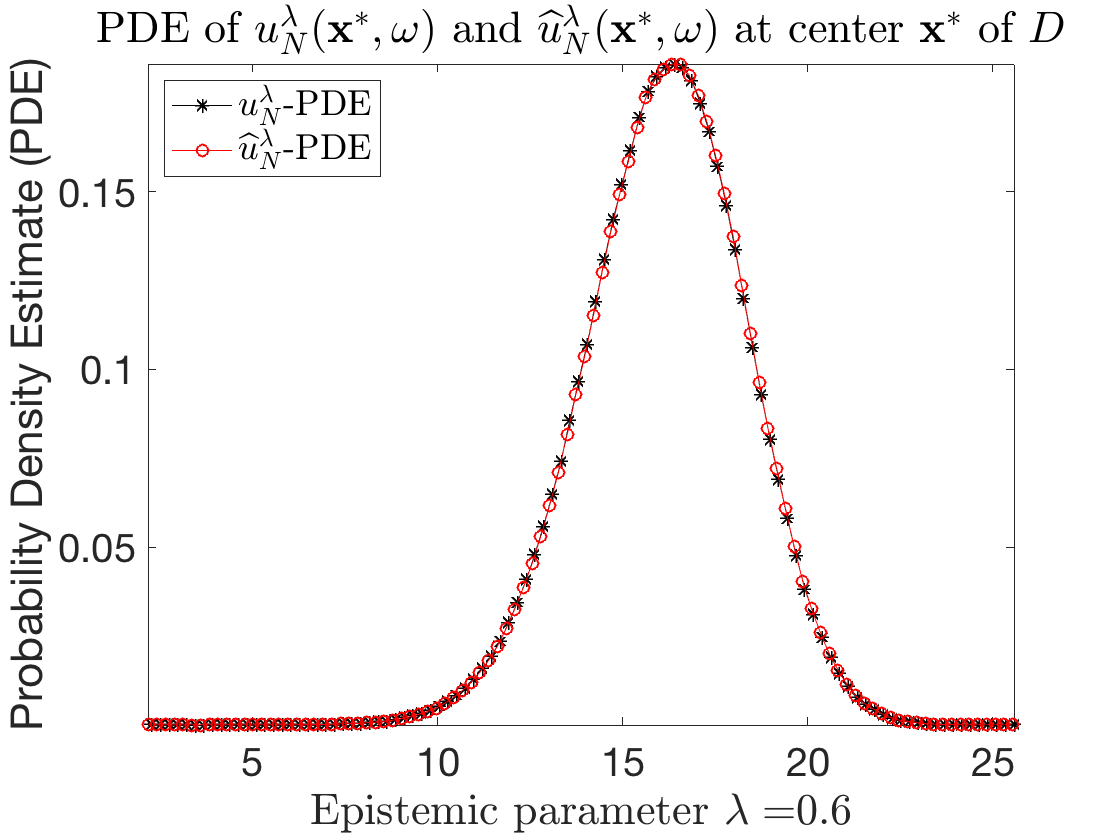}\\
    \parbox{4cm}{\centering
      (d) $\tau = 0.4$}
    &
    \parbox{4cm}{\centering
      (e) $\tau = 0.5$}
    &
    \parbox{4cm}{\centering
      (f) $\tau = 0.6$} \\\includegraphics[width=4cm]{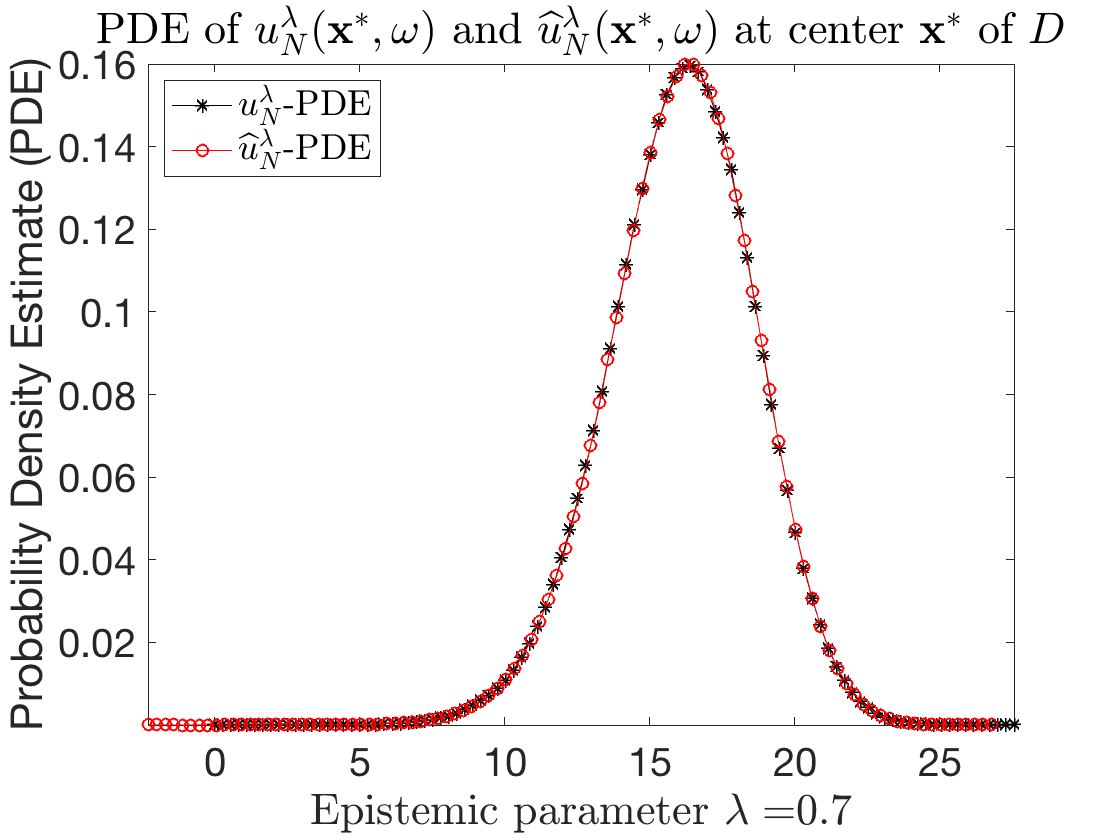}
    &
    \includegraphics[width=4cm]{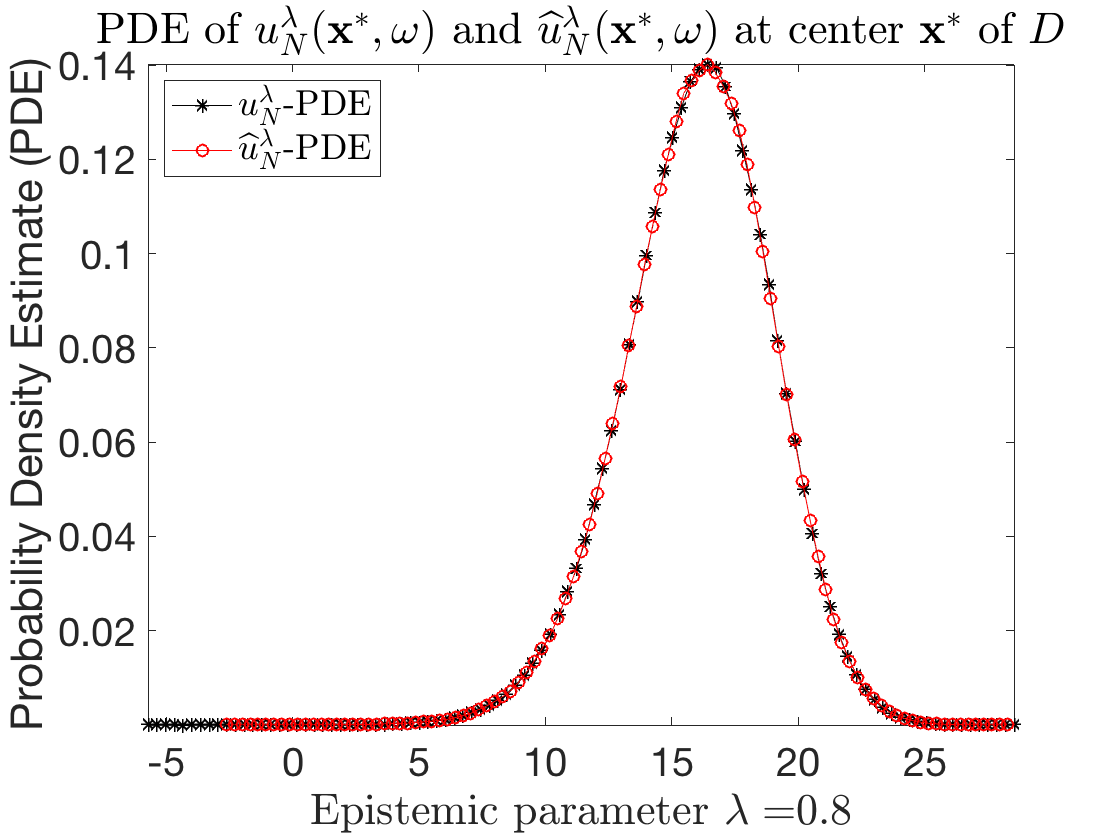}
    &
    \includegraphics[width=4cm]{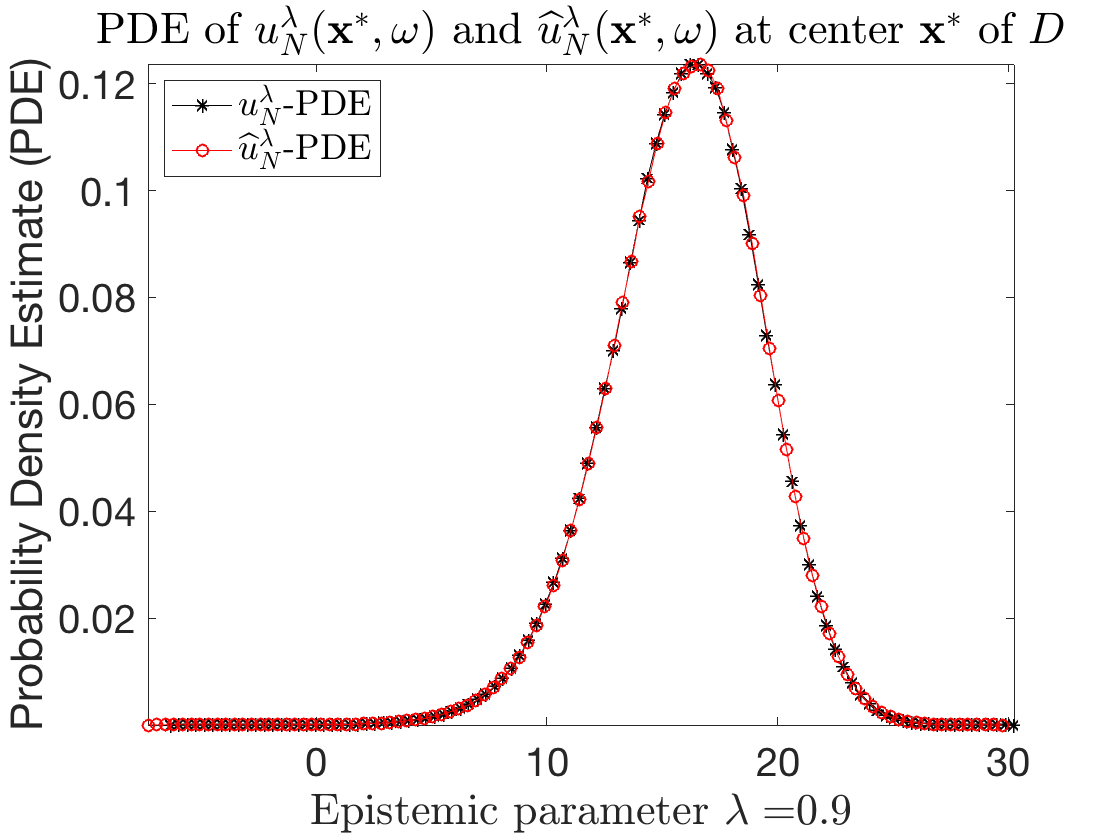}\\
    \parbox{4cm}{\centering
      (g) $\tau = 0.7$}
    &
    \parbox{4cm}{\centering
      (h) $\tau = 0.8$}
    &
    \parbox{4cm}{\centering
      (i) $\tau = 0.9$} \\
  \end{tabular}
  \caption{\label{fig:PDE-d10}
  {\bf Ten-dimensional stochastic diffusion model}:
  Comparisons of  probability density estimates (PDE) of the gPC solution $u_N^\tau(\boldsymbol{x}^*, \omega)$
  and corresponding post-processed surrogate $\widehat{u}_N^\tau(\boldsymbol{x}^*, \omega)$ at the center $\boldsymbol{x}^*$ of $D$,
  and for various epistemic parameter values.} 
\end{figure}

\begin{figure}[]
  \centering
  \begin{tabular}{ccc}
    \includegraphics[width=4cm]{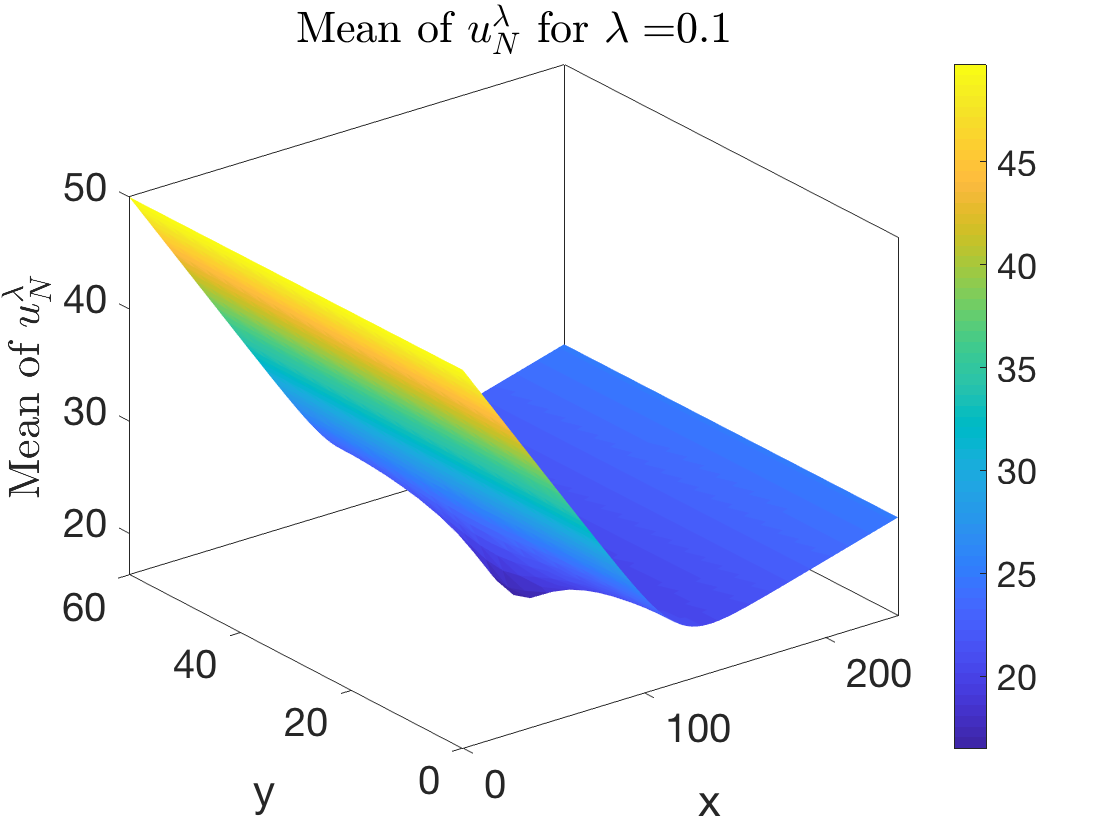}
    &
    \includegraphics[width=4cm]{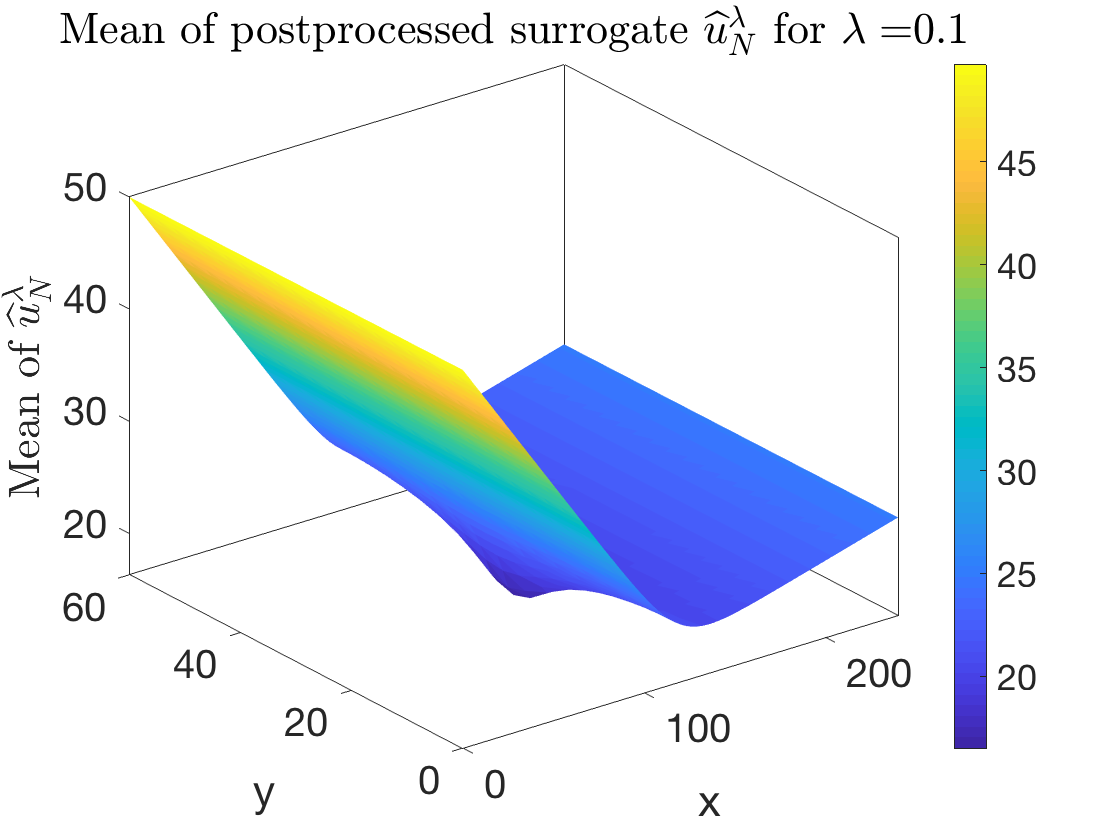}
    &
    \includegraphics[width=4cm]{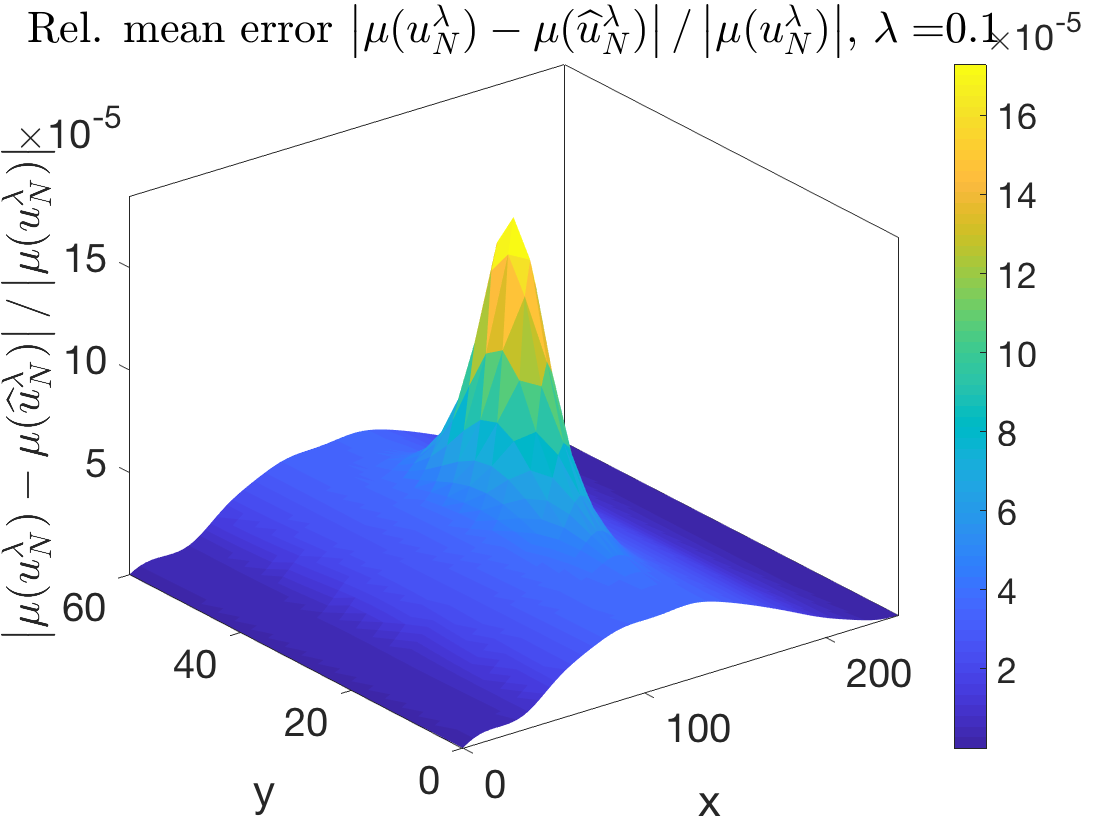}\\
    \parbox{4cm}{\centering
      (a) $\tau = 0.1$, \\ gPC mean solution $u_N^\tau$}
    &
    \parbox{4cm}{\centering
      (b) $\tau = 0.1$, \\ post-processed surrogate  mean solution  $\widehat{u}_N^\tau$}
    &
    \parbox{4cm}{\centering
      (c) $\tau = 0.1$, \\ relative error of mean solutions} \      
       \end{tabular}
  \caption{\label{fig:soln-d10}
  {\bf Ten-dimensional stochastic diffusion model}:
  Comparisons of  the mean gPC solution $u_N^\tau$
  and the mean of corresponding post-processed surrogate $\widehat{u}_N^\tau$ 
  for $\tau = 0.1$.} 
\end{figure}

For the ten-dimensional stochastic diffusion model, 
results in Table~\ref{eq:table-epistemic-err-d10} demonstrate the  accuracy  in
mean and variance of the post-processed surrogate solution $\widehat{u}_N^\tau$ when compared
to the gPC reference solution $u_N^\tau$, for $\tau = 0.1, 0.2, \dots, 0.9$. The maximum
relative error was computed by taking the maximum of the relative errors at the spatial grid points.
In Figure~\ref{fig:PDE-d10} we demonstrate the accuracy of the post-processing
approach by comparing the
probability density estimates for
$\widehat{u}_N^\tau(\boldsymbol{x}^*, \cdot)$ and
$u_N^\tau(\boldsymbol{x}^*, \cdot)$
for $\boldsymbol{x}^*$ at the centre of $D$.
A  sample profile in $D$, for the
case $\tau = 0.1$, of the
ten-dimensional mean solution $u_N^\tau(\cdot, \omega)$,
the counterpart post-processed surrogate $\widehat{u}_N^\tau(\cdot, \omega)$, 
and the associated relative error, the ten-dimensional EUQ stochastic model  are given in 
Figure~\ref{fig:soln-d10}. 

\begin{figure}
  \centering
  \includegraphics[width=10cm]{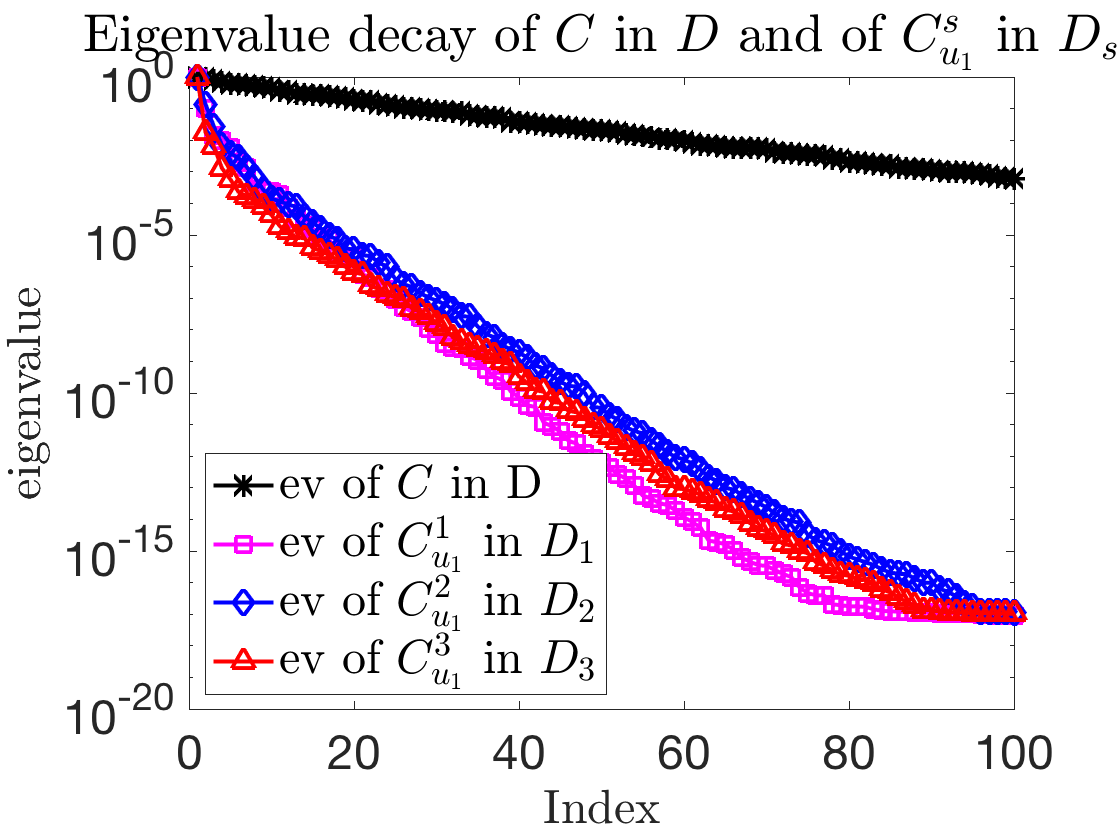}
  \caption{\label{fig:eigenvalue-d100}
    Decay of the eigenvalues of
    the given covariance matrix $C(\x,\y)$  at grid points in $D \times D$
    and the coarse solution $u_1$ based covariance matrix $C^s_{u_1}(\x,\y)$ at
    grid points in the sub-domains $D_s \times D_s$, for $s = 1, 2, 3$, where
    $D = \cup_{s=1}^8 D_s$.
  }
\end{figure}

The decay of the eigenvalues of $C$  in $D$ shown in 
Figure~\ref{fig:eigenvalue-d100} suggests that it may be appropriate
to increase the the truncation parameter for the 
EUQ model to $d =100$. We note that it is not feasible to solve  the
gPC model in 100-dimensions, with polynomial degree $N = 3$,
because the PDE would need to be solved for billons of sparse grid
realization points.
However,
using the gPC-DD approach it is feasible.
In particular, we
first solve for $u_1$ in 100-dimensions
(with $N =1$,  and sparse-grid level $2$), which requires 
solving the PDE only $201$  times, corresponding to the sparse-grid points
in $100$-dimensions.
We then use $u_1$
in conjunction with the domain-decomposition method with $S  = 8$ subdomains
of $D = (0,240) \times (0,60)$,
in a $4 \times 2$ grid.

For $s=1, \dots, 8$, the eigenvalues of $C_{u_1}^s$ in $D_s$ decay substantially faster, than those
of $C$ in $D$. This is demonstrated (for $s = 1, 2, 3$) in Figure~\ref{fig:eigenvalue-d100}. Hence, 
after solving for $u_1$ in the 100-dimensional stochastic space we
solve the sub-domain problems
using $r = 3$ stochastic dimensional spaces.
The transfer between the sub-domain stochastic variables and the
full-domain stochastic variables are computed using the  representation in~\eqref{eq:eta-to-xi}.
For $0  < \tau < 1$ we carry out the post-processing surrogate approach
in Section~\ref{sec:epistemic-dd}
for each sub-domain model, and compute the corresponding  mean and variance
efficiently using~\eqref{eq:mean-dd-surrogate}.

\begin{figure}
  \centering
  \begin{tabular}{ccc}
    \includegraphics[width=4cm]{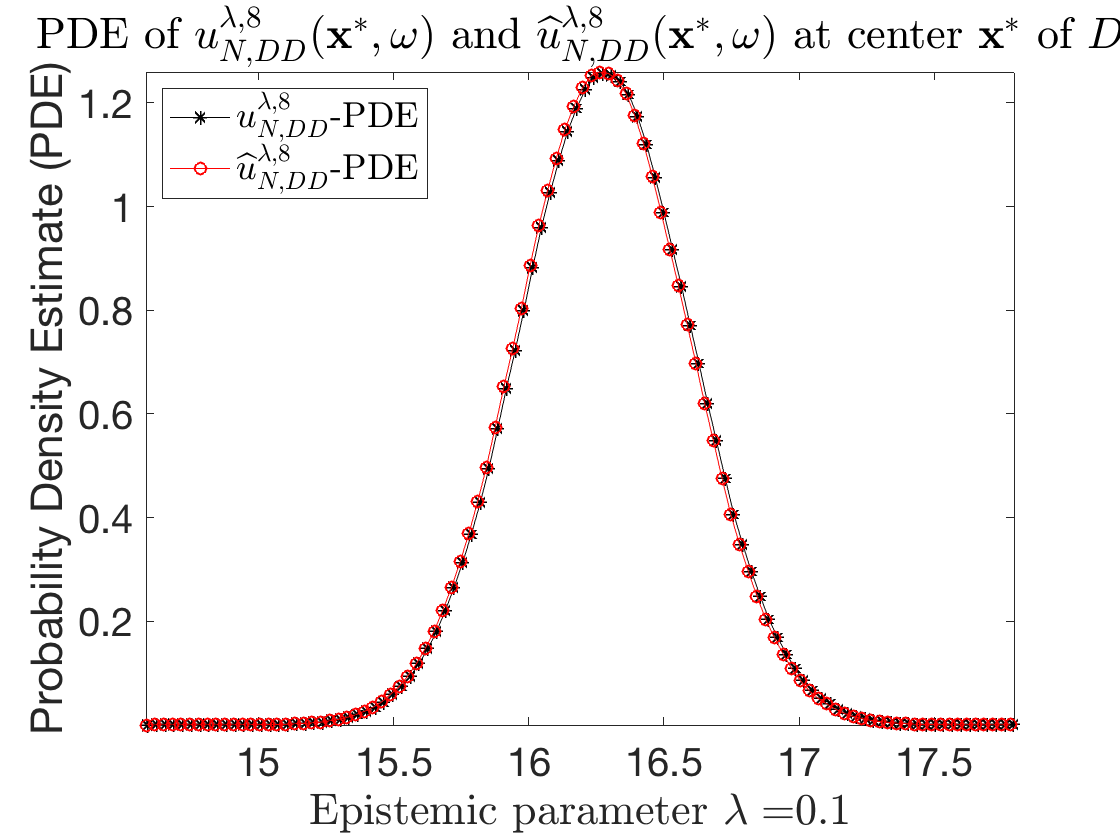}
    &
    \includegraphics[width=4cm]{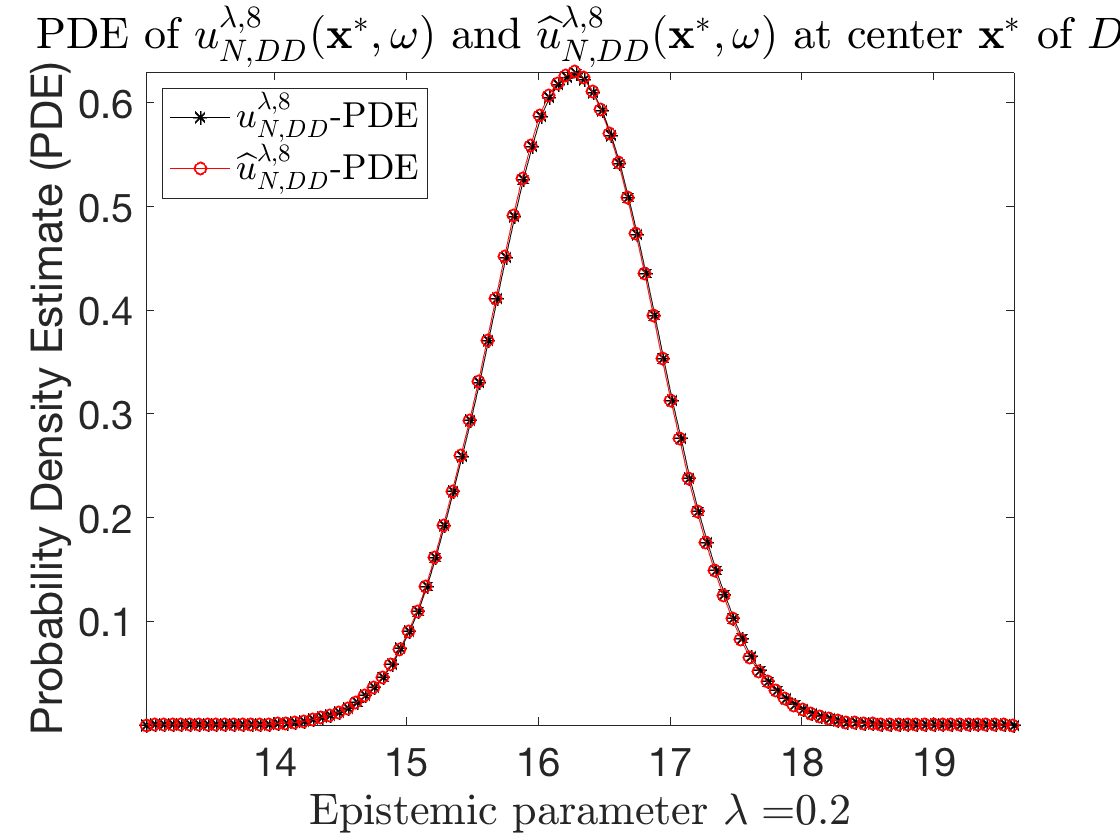}
    &
    \includegraphics[width=4cm]{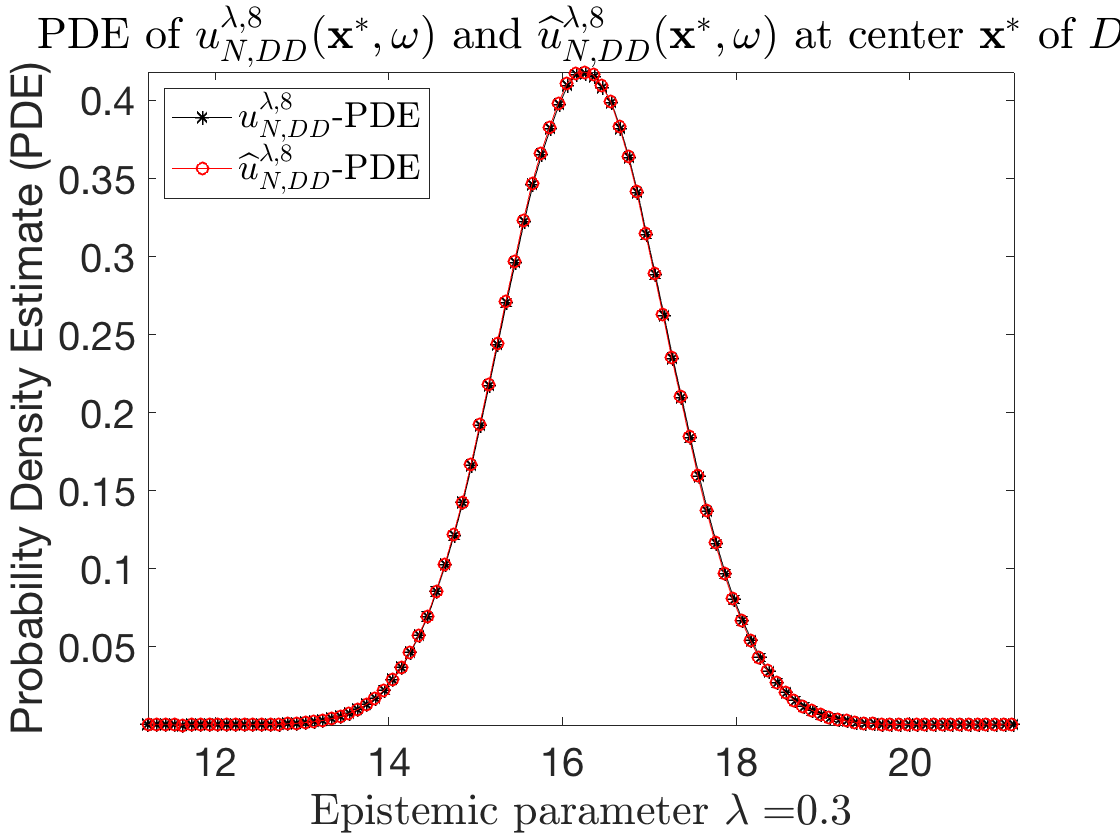}\\
    \parbox{4cm}{\centering
      (a) $\tau = 0.1$}
    &
    \parbox{4cm}{\centering
      (b) $\tau = 0.2$}
    &
    \parbox{4cm}{\centering
      (c) $\tau = 0.3$} \\
      \includegraphics[width=4cm]{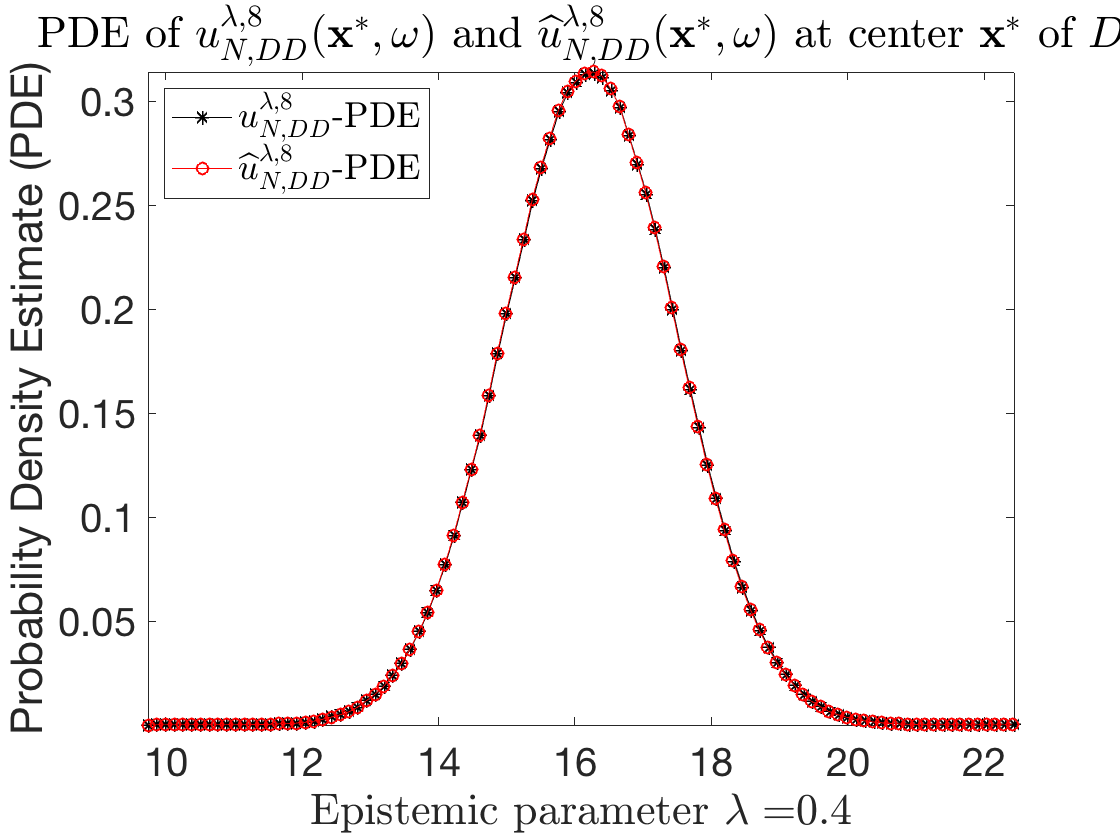}
    &
    \includegraphics[width=4cm]{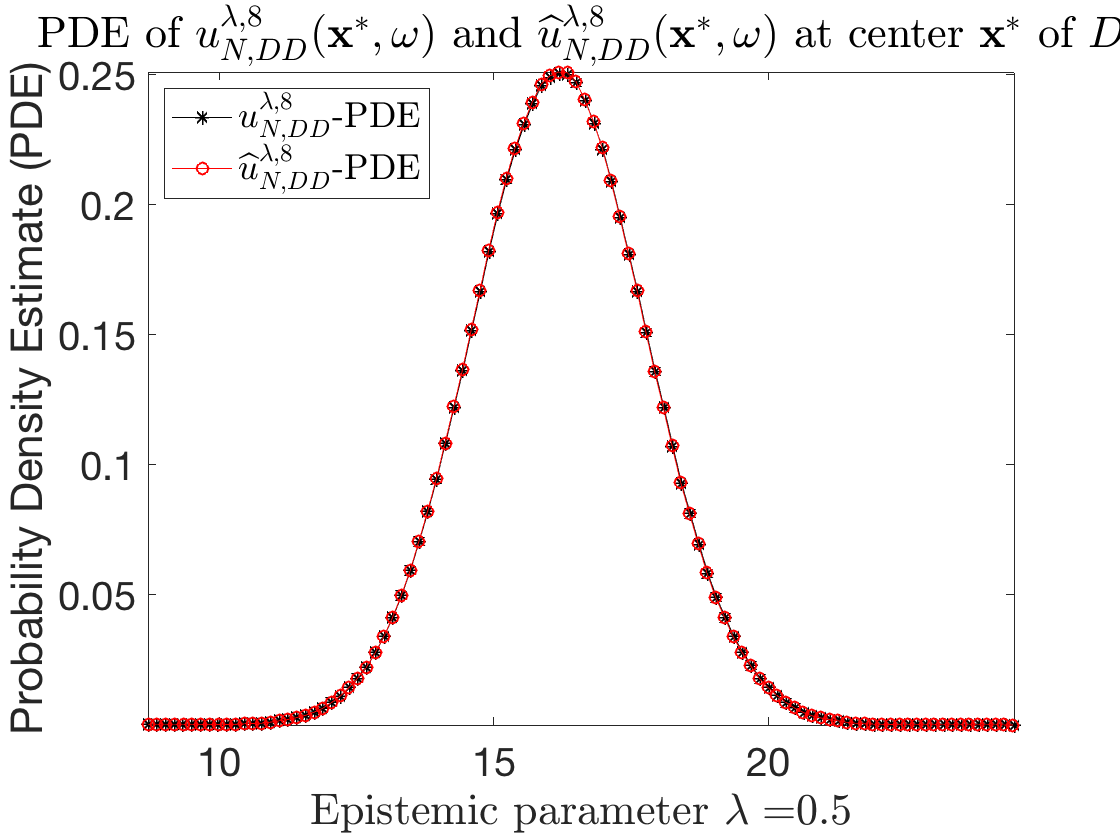}
    &
    \includegraphics[width=4cm]{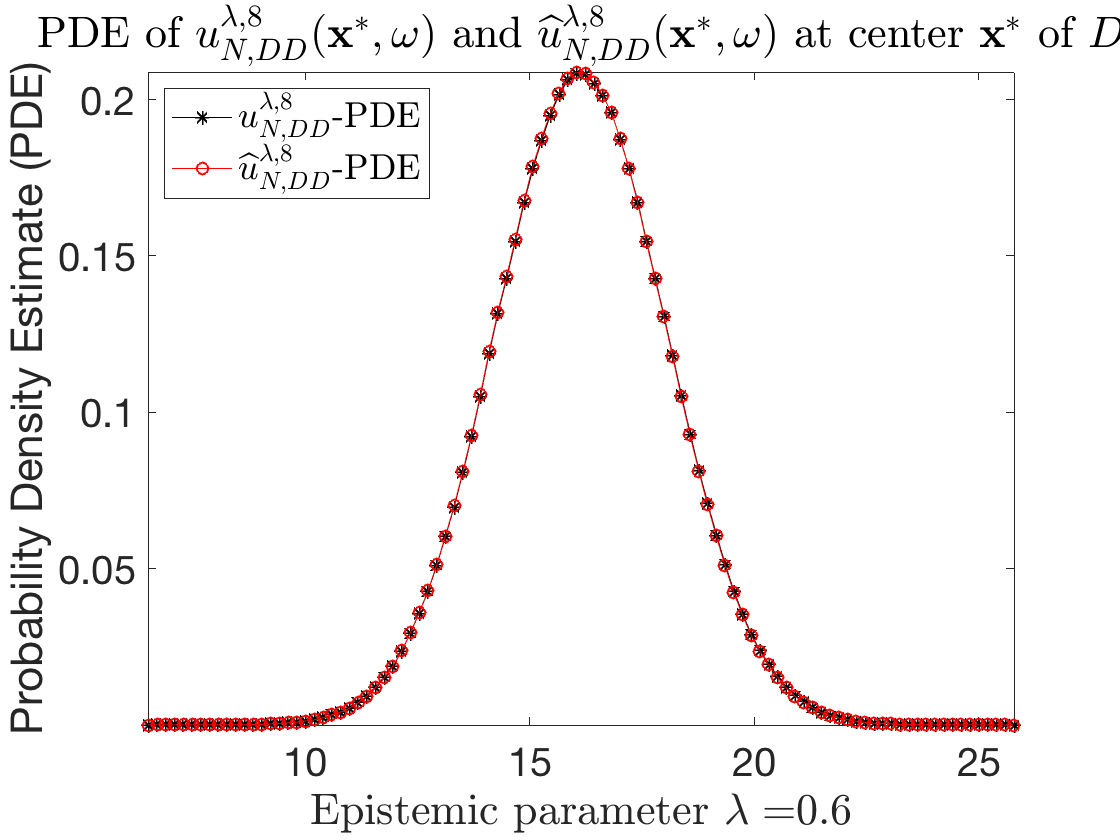}\\
    \parbox{4cm}{\centering
      (d) $\tau = 0.4$}
    &
    \parbox{4cm}{\centering
      (e) $\tau = 0.5$}
    &
    \parbox{4cm}{\centering
      (f) $\tau = 0.6$} \\\includegraphics[width=4cm]{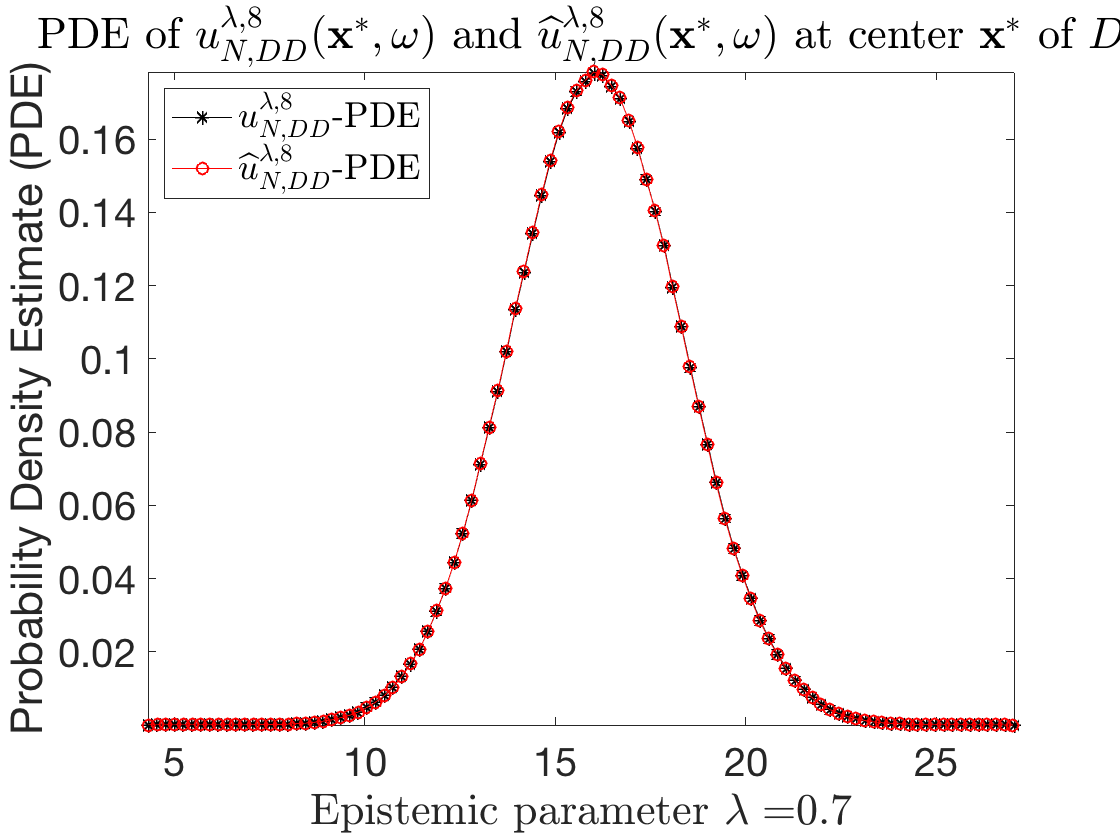}
    &
    \includegraphics[width=4cm]{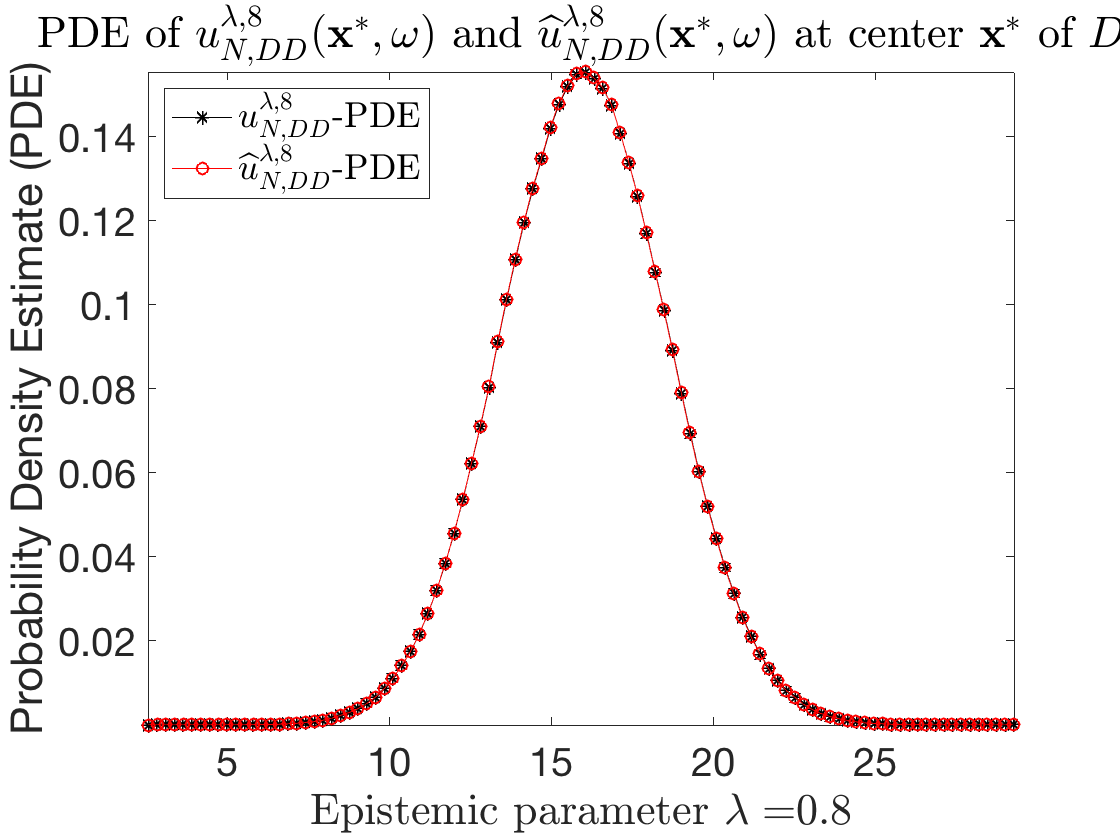}
    &
    \includegraphics[width=4cm]{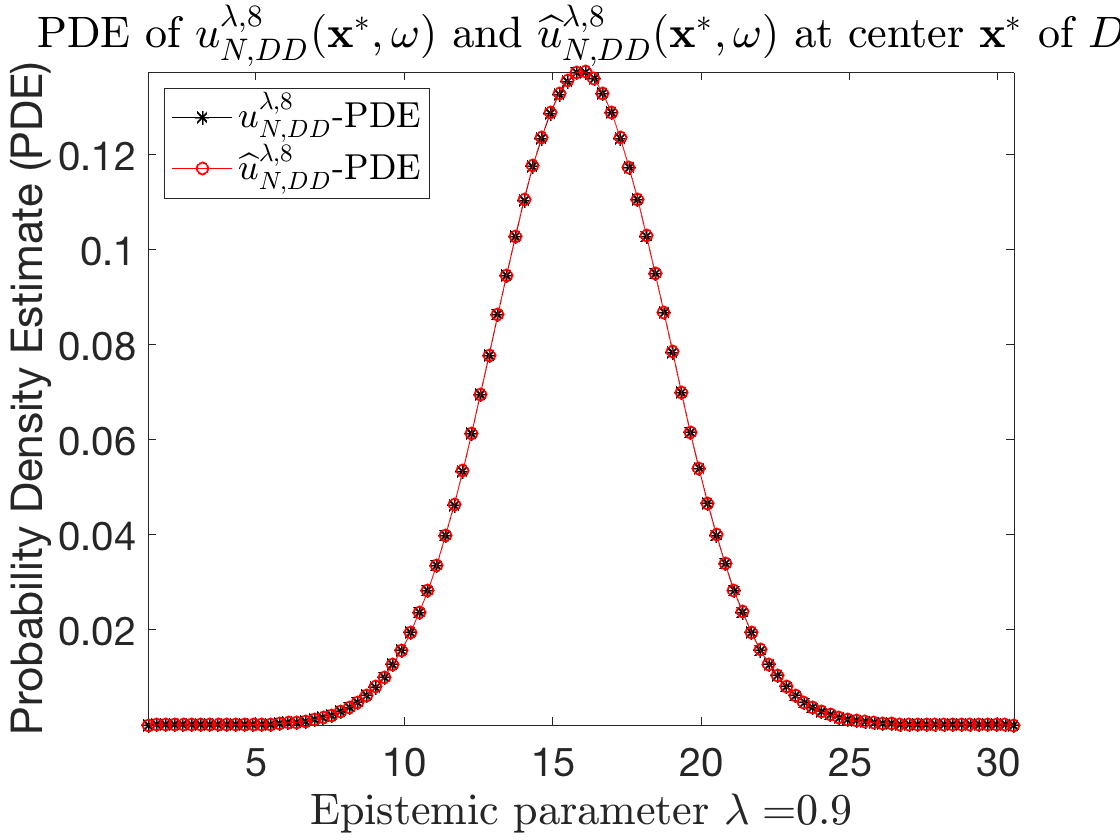}\\
    \parbox{4cm}{\centering
      (g) $\tau = 0.7$}
    &
    \parbox{4cm}{\centering
      (h) $\tau = 0.8$}
    &
    \parbox{4cm}{\centering
      (i) $\tau = 0.9$} \\
  \end{tabular}
  \caption{\label{fig:PDE-d100}
  {\bf Hundred-dimensional stochastic diffusion model}:
  Comparisons of  probability density estimates (PDE) of the gPC-DD solution $u_{N,DD}^{\tau,8}(\boldsymbol{x}^*, \omega)$
  and corresponding post-processed surrogate $\widehat{u}_{N,DD}^{\tau,8}(\boldsymbol{x}^*, \omega)$
  at the center $\boldsymbol{x}^*$ of $D$,
  and for various epistemic parameter values.} 
\end{figure}

\begin{figure}
  \centering
  \begin{tabular}{ccc}
    \includegraphics[width=4cm]{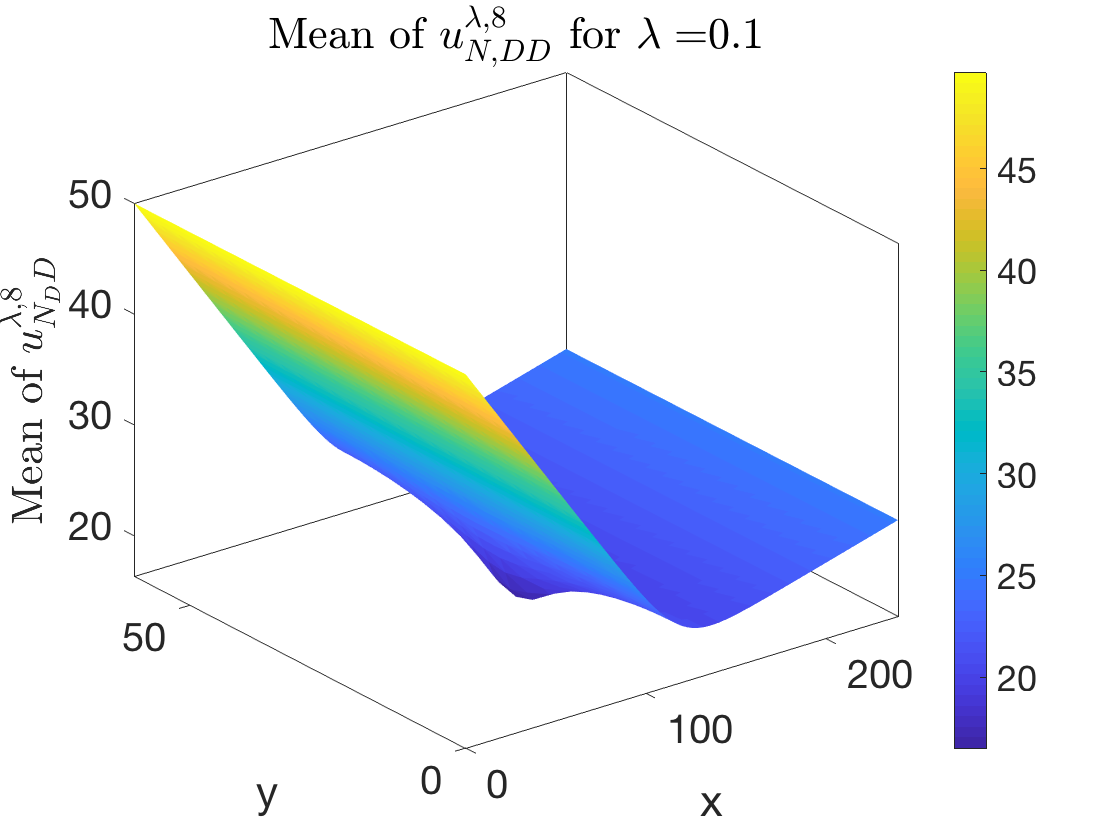}
    &
    \includegraphics[width=4cm]{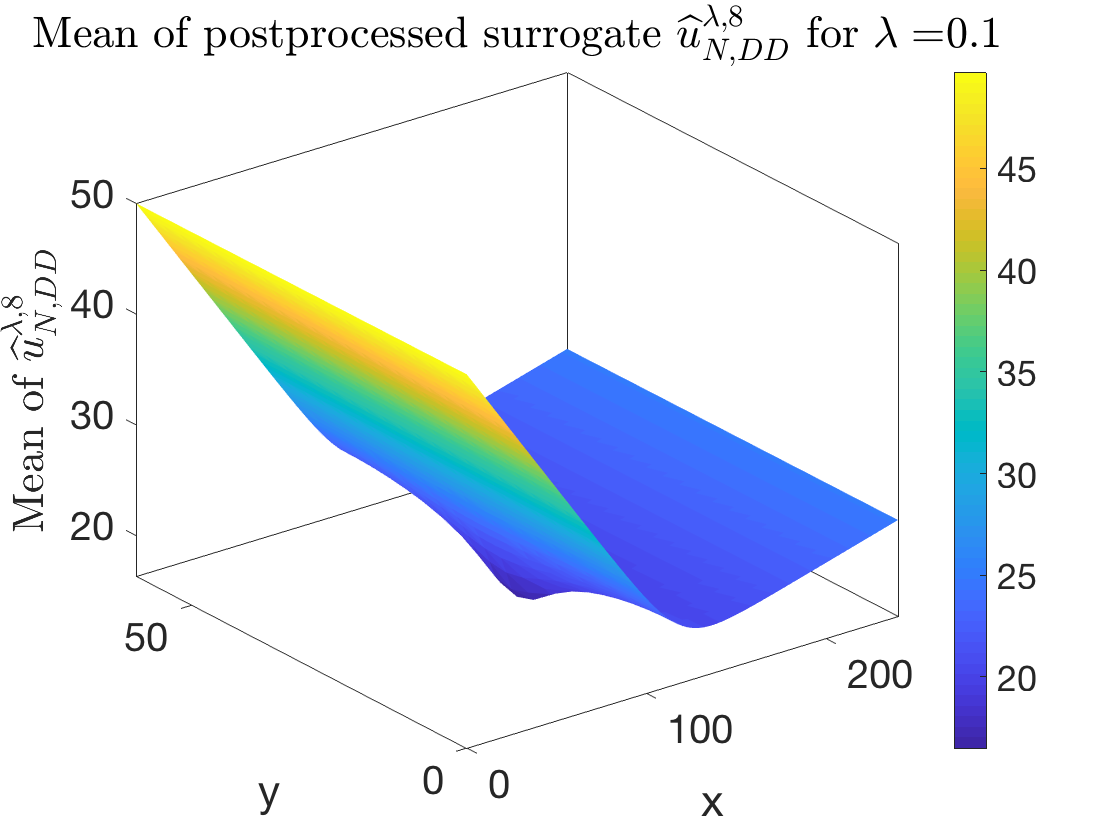}
    &
    \includegraphics[width=4cm]{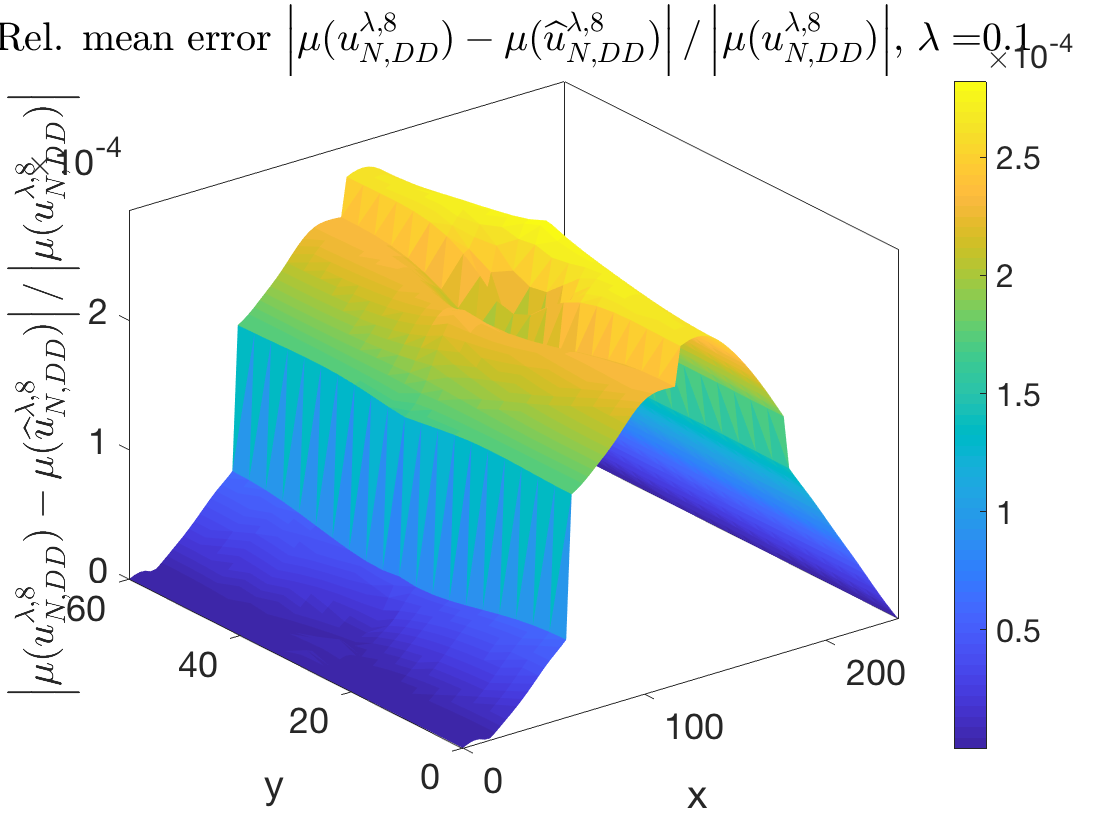}\\
    \parbox{4cm}{\centering
      (a) $\tau = 0.1$, \\ gPC-DD mean solution $u_{N,DD}^{\tau,8}$}
    &
    \parbox{4cm}{\centering
      (b) $\tau = 0.1$, \\ post-processed surrogate  mean solution  $\widehat{u}_{N,DD}^{\tau,8}$}
    &
    \parbox{4cm}{\centering
      (c) $\tau = 0.1$, \\ relative error of mean solutions} \      
       \end{tabular}
  \caption{\label{fig:soln-d100}
  {\bf Hundred-dimensional stochastic diffusion model}:
  Comparisons of  the mean gPC-DD solution $u_{N,DD}^{\tau,8}$
  and the mean of corresponding post-processed surrogate 
  $\widehat{u}_{N,DD}^{\tau,8}$
  for $\tau = 0.1$.} 
\end{figure}

For the hundred-dimensional stochastic diffusion models, 
results in Table~\ref{eq:table-epistemic-err-d100} demonstrate the  accuracy  in
mean and variance of the post-processed surrogate solution $\widehat{u}_{N,DD}^{\tau,8}$ when compared
to the gPC-DD reference solution $u_{N,DD}^{\tau,8}$, for $\tau = 0.1, 0.2, \dots, 0.9$. The maximum
relative error was computed by taking the maximum of the relative errors at the spatial grid points. 
In Figure~\ref{fig:PDE-d100} we demonstrate the accuracy of the post-processing
approach by comparing the
probability density estimates for
$\widehat{u}_{N,DD}^{\tau,8}(\boldsymbol{x}^*, \cdot)$ and
${u}_{N,DD}^{\tau,8}(\boldsymbol{x}^*, \cdot)$
for $\boldsymbol{x}^*$ at the centre of $D$. 
In  Figure~\ref{fig:soln-d100}
we visualize the mean $u_{N,DD}^{\tau,8}(\cdot, \omega)$ and
the post-processed surrogate approximation
$\widehat{u}_{N,DD}^{\tau,8}(\cdot, \omega)$,
for the hundred-dimensional EUQ stochastic model with $\tau = 0.1$.

\begin{table}[!ht]
  \centering
  \begin{tabular}{cccc}
    \hline
    $\tau$ & max. rel.~mean error & max. rel.~variance error\\
    \hline
    0.9 &	6.9082e-06  & 1.6753e-04 \\
    0.8 &	2.4649e-05   & 3.0185e-04  \\
    0.7 &	4.9307e-05 &  4.0997e-04 \\
    0.6 &	7.7467e-05  & 4.9846e-04 \\
    0.5  &	1.0619e-04  & 5.7096e-04\\
    0.4 &	1.3302e-04 &  7.1537e-04 \\
    0.3 &	1.5595e-04  & 9.0909e-04 \\
    0.2 &	1.7344e-04  & 1.0596e-03 \\
    0.1 &	 1.8436e-04  & 1.1548e-03 \\

    \hline
  \end{tabular}
  \caption{\label{eq:table-epistemic-err-d10}
    {\bf Ten-dimensional stochastic diffusion model}:
    Maximum norm relative error of the approximate mean and variance computed using the
    gPC approach and the post-processed   surrogate method.} 
\end{table}

  \clearpage

 \clearpage

\begin{table}[!ht]
  \centering
  \begin{tabular}{cccc}
    \hline
    $\tau$ & max. rel.~mean error & max. rel.~variance error\\
    \hline
    0.9 &	8.8316e-05  &  1.2265e-03  \\
    0.8 &	1.5066e-04   & 2.2016e-03  \\
    0.7 &	1.9168e-04 &  2.9742e-03 \\
    0.6 &	2.2377e-04  & 3.8097e-03  \\
    0.5  &	2.4742e-04 &  4.5387e-03\\
    0.4 &	 2.6353e-04 &  5.1387e-03 \\
    0.3 &	 2.7460e-04  & 5.6057e-03 \\
    0.2 &	2.8140e-04  & 5.9388e-03 \\
    0.1 &	2.8503e-04   &  6.1383e-03 \\

    \hline
  \end{tabular}
  \caption{\label{eq:table-epistemic-err-d100}
    {\bf Hundred-dimensional stochastic diffusion model}:
    Maximum norm relative error of the approximate mean and variance computed using the
    DD-gPC approach and the post-processed   surrogate method.} 
\end{table}

\newpage

\clearpage

\section*{Acknowledgments}
This research was supported by the U.S. Department of Energy (DOE),  Office of Advanced 
Scientific Computing Research.  Pacific Northwest National Laboratory is operated by Battelle for the DOE under Contract DE-AC05-76RL01830.

\bibliography{references}

\end{document}